\documentclass{article}

\usepackage{amsfonts}
\usepackage{amsmath}
\usepackage{amssymb}
\usepackage{graphicx}
\usepackage{enumitem}

\newcommand{\subjclass}[1]{\par\noindent{\bf MSC 2020: }#1\par}
\newcommand{\keywords}[1]{\par\noindent{\bf Keywords: }#1\par}

\providecommand{\U}[1]{\protect\rule{.1in}{.1in}}

\newtheorem{theorem}{Theorem}

\newtheorem{corollary}[theorem]{Corollary}

\newtheorem{definition}[theorem]{Definition}

\newtheorem{lemma}[theorem]{Lemma}

\newtheorem{remark}[theorem]{Remark}

\newenvironment{proof}[1][Proof]{\noindent\textbf{#1.} }{\ \rule{0.5em}{0.5em}}

\title{Robustness of dynamically gradient multivalued dynamical systems}

\author{
Rub\'en Caballero\thanks{Centro de Investigación Operativa, Universidad Miguel Hernández de Elche, 03202-Elche, Alicante, Spain. 
Email: \texttt{ruben.caballero@umh.es}} 
\and
Alexandre N. Carvalho\thanks{Instituto de Ciências Matemáticas e de Computação, Universidade de São Paulo, Caixa Postal 668, 13560-970 São Carlos, SP, Brazil. 
Email: \texttt{andcarva@icmc.usp.br}} 
\and
Pedro Mar\'in-Rubio\thanks{Departamento de Ecuaciones Diferenciales y Análisis Numérico, Universidad de Sevilla, 41012-Sevilla, Spain.
Email: \texttt{pmr@us.es}}
\and
Jos\'e Valero\thanks{Centro de Investigación Operativa, Universidad Miguel Hernández de Elche, 03202-Elche, Alicante, Spain.
Email: \texttt{jvalero@umh.es}}
\and
\thanks{The authors of this work have been partially funded by the following projects:
R. Caballero is a fellow of Programa de FPU del Ministerio de Educación, Cultura y Deporte, reference FPU15/03080;
P. Mar\'in-Rubio was supported by projects PHB2010-0002-PC (Ministerio de Educación-DGPU), MTM2015-63723-P (MINECO/FEDER, EU), and P12-FQM-1492 (Junta de Andalucía);
and J. Valero by projects MTM2015-63723-P, MTM2016-74921-P (MINECO/FEDER, EU) and P12-FQM-1492 (Junta de Andalucía).}
}

\date{}

\begin{document}

\maketitle


{\footnotesize \centerline{To Professor Valery Melnik, in Memoriam}}

\begin{abstract}
In this paper we study the robustness of dynamically gradient multivalued semiflows. As an application, we describe the dynamical properties of a family of Chafee-Infante problems approximating a differential inclusion studied in \cite{arrietavalero}, proving that the weak solutions of these problems generate a dynamically gradient multivalued semiflow with respect to suitable Morse sets.
\end{abstract}

\bigskip

\section{Introduction}

One of the main goals of the theory of dynamical systems is to characterize
the structure of global attractors. It is possible to find a wide literature
about this problem for semigroups; however, it has been recently when new
results in this direction for multivalued dynamical systems have been proved
\cite{arrietavalero}, \cite{kapustyankasyanov},
\cite{kapustyankasyanovvalero2015}.

In this sense, the theory of Morse decomposition plays an important role. In
fact, the existence of a Lyapunov function, the property of being a
dynamically gradient semiflow and the existence of a Morse decomposition are
shown to be equivalent for multivalued dynamical systems in \cite{costavalero}.

In this work we show under suitable assumptions that a dynamically gradient
multivalued semiflow is stable under perturbations, that is, the family of
perturbed multivalued semiflows remains dynamically gradient.

For a fixed dynamically gradient multivalued semiflow with a global attractor
we also analyze the rearrangement of a pairwise disjoint finite family of
isolated weakly invariant sets, included in the attractor, in such a way that
the dynamically gradient property is satisfied in the stronger sense of
\cite{li}.

These results extend previous ones in the single-valued framework in
\cite{carvalholanga, aragaocosta2011, aragaocosta2013} to the case where
uniqueness of solution does not hold. Additionally, it is worth saying that
the m-semiflows here are not supposed to be general dynamical systems as in
\cite{li}, where a robustness theorem for Morse decompositions of multivalued
dynamical systems is also proved under a suitable continuity assumption.

We also apply this general robustness theorem in order to show that a family
of Chafee-Infante problems approximating a differential inclusion is
dynamically gradient it it is close enough to the original problem.

This paper is organized as follows.

Firstly, we introduce in Section \ref{preliminaries} basic concepts and
properties related to fixed points, complete trajectories and global
attractors. In this way we are able to present in Section \ref{Robustness} the
main result about robustness of dynamically gradient multivalued semiflows.
Further, in Section \ref{equivalent} we prove a theorem which allows us to
reorder the family of weakly invariants sets, thus establishing an equivalent
definition of dynamically gradient families.

Afterwards, we consider a Chafee-Infante problem in Section \ref{application},
where the equivalence of weak and strong solutions is established. Once the
set of fixed points is analyzed, we consider a family of Chafee-Infante
equations, approximating the differential inclusion tackled in
\cite{arrietavalero}. We check that this family of Chafee-Infante equations
verifies the hypotheses of the robustness theorem in order to obtain,
therefore, that the multivalued semiflows generated by the solutions of the
approximating problems are dynamically gradient if this family is closed
enough to the original one.

\section{Preliminaries}

\label{preliminaries} Consider a metric space $(X,d)$ and a family of
functions $\mathcal{R} \subset\mathcal{C} (\mathbb{R}_{+}; X) $. Denote by
$P(X)$ the class of nonempty subsets of $X$. Then, define the multivalued map
$G: \mathbb{R}_{+} \times X \rightarrow P(X)$ associated with the family
$\mathcal{R}$ as follows
\begin{equation}
\label{defG}G(t,u_{0})=\{u(t) : u(\cdot) \in\mathcal{R}, u(0)=u_{0} \}.
\end{equation}

In this abstract setting, the multivalued map $G$ is expected to satisfy some
properties that fit in the framework of multivalued dynamical systems. The
first concept is given now, although a more axiomatic construction will be
provided below.

\begin{definition}
Let $(X,d)$ be a metric space. A multivalued map $G:\mathbb{R}_{+} \times X
\rightarrow P(X)$ is a multivalued semiflow (or m-semiflow) if $G(0,x)=x$ for
all $x \in X$ and $G(t+s,x) \subset G(t,G(s,x))$ for all $t,s\geq0$ and $x \in
X$. \newline If the above is not only an inclusion, but an equality, it is
said that the m-semiflow is strict.
\end{definition}

Once a multivalued dynamical system is defined, we recall the concepts of
invariance and global attractor, with evident differences with respect to the
single-valued case.

\begin{definition}
A map $\gamma: \mathbb{R }\rightarrow X $ is called a complete trajectory of
$\mathcal{R}$ (resp. of G) if $\gamma(\cdot+h)\mid_{[0,\infty)} \in
\mathcal{R}$ for all $h \in\mathbb{R}$ (resp. if $\gamma(t+s) \in G(t,
\gamma(s))$ for all $s \in\mathbb{R}$ and $t \geq0)$.

A point $z \in X$ is a fixed point of $\mathcal{R}$ (resp. of G) if
$\varphi(\cdot)\equiv z \in\mathcal{R}$ (resp. $z \in G(t,z)$ for all $t
\geq0$).
\end{definition}

\begin{definition}
Given an m-semiflow $G$ on a metric space $(X,d)$ a set $B \subset X$ is said
to be negatively invariant if $B \subset G(t,B)$ for all $t\geq0$, and
strictly invariant (or, simply, invariant) if the above relation is not only
an inclusion but an equality.

The set $B$ is said to be weakly invariant if for any $x \in B$ there exists a
complete trajectory $\gamma$ of $\mathcal{R}$ contained in $B$ such that
$\gamma(0) =x $. We observe that weakly invariance implies negatively invariance.

A set $\mathcal{A} \subset X$ is called a global attractor for an m-semiflow
if it is negatively semi-invariant, i.e., $\mathcal{A} \subset G(t,\mathcal{A}%
) $ for all $t \geq0$, and it attracts all attainable sets through the
m-semiflow starting in bounded subsets, i.e., $dist_{X}(G(t,B),\mathcal{A})
\rightarrow0$ as $t \rightarrow\infty$, where $dist_{X}(A,B)= \sup_{a \in A}
\inf_{b \in B} d(a,b)$.
\end{definition}

\begin{remark}
A global attractor for an m-semiflow does not have to be unique, nor a bounded
set. However, if a global attractor is bounded and closed, it is minimal among
all closed sets that attract bounded sets \cite{melnikvalero}. In particular,
a bounded and closed global attractor is unique.
\end{remark}

In order to obtain a detailed characterization of the internal structure of a
global attractor, we introduce an axiomatic set of properties on the set
$\mathcal{R}$ (see \cite{ball} and \cite{kapustyankasyanov}).

The set of axiomatic properties that we will deal with is the following.

\begin{itemize}
\item[(K1)] For any $x \in X$ there exists at least one element $\varphi
\in\mathcal{R}$ such that $\varphi(0) =x$.

\item[(K2)] $\varphi_{\tau} (\cdot) := \varphi(\cdot+\tau) \in\mathcal{R}$ for
any $\tau\geq0$ and $\varphi\in\mathcal{R}$ (translation property).

\item[(K3)] Let $\varphi_{1}, \varphi_{2} \in\mathcal{R}$ be such that
$\varphi_{2}(0)=\varphi_{1}(s)$ for some $s>0$. Then, the function $\varphi$
defined by
\[
\varphi(t) = \left\lbrace
\begin{array}
[c]{lcc}%
\varphi_{1}(t) \quad0 \leq t \leq s &  & \\
\varphi_{2}(t-s) \quad s \leq t, &  & 
\end{array}
\right.
\]
belongs to $\mathcal{R}$ (concatenation property).

\item[(K4)] For any sequence $\{\varphi^{n}\} \subset\mathcal{R}$ such that
$\varphi^{n}(0) \rightarrow x_{0}$ in X, there exist a subsequence
$\{\varphi^{n_{k}}\}$ and $\varphi\in\mathcal{R}$ such that $\varphi^{n_{k}%
}(t) \rightarrow\varphi(t)$ for all $t \geq0$.
\end{itemize}

It is immediate to observe \cite[Proposition 2]{caraballorubio} or \cite[Lemma
9]{kapustyanpankov} that $\mathcal{R}$ fulfilling (K1) and (K2) gives rise to
an m-semiflow $G$ through (\ref{defG}), and if besides (K3) holds, then this
m-semiflow is strict. In such a case, a global bounded attractor, supposing
that it exists, is strictly invariant \cite[Remark 8]{melnikvalero}.

Several properties concerning fixed points, complete trajectories and global
attractors are summarized in the following results \cite{kapustyankasyanov}.

\begin{lemma}
Let (K1)-(K2) be satisfied. Then every fixed point (resp. complete trajectory)
of $\mathcal{R}$ is also a fixed point (resp. complete trajectory) of G.

If $\mathcal{R}$ fulfills (K1)-(K4), then the fixed points of $\mathcal{R}$
and G coincide. Besides, a map $\gamma:\mathbb{R} \rightarrow X$ is a complete
trajectory of $\mathcal{R}$ if and only if it is continuous and a complete
trajectory of G.
\end{lemma}

The standard well-known result in the single-valued case for describing the
attractor as the union of bounded complete trajectories reads in the
multivalued case as follows.

\begin{theorem}
\label{structureattractor} Consider $\mathcal{R}$ satisfying (K1) and (K2),
and either (K3) or (K4). Assume also that G possesses a compact global
attractor $\mathcal{A}$. Then
\begin{equation}
\label{attractor}\mathcal{A}= \{ \gamma(0) : \gamma\in\mathbb{K} \}=
\displaystyle \cup_{t \in\mathbb{R}} \{ \gamma(t) : \gamma\in\mathbb{K} \},
\end{equation}
where $\mathbb{K}$ denotes the set of all bounded complete trajectories in
$\mathcal{R}$.
\end{theorem}

Now we recall the definitions of some important sets in the literature of the
dynamical systems. Let $B\subset X$ and let $\varphi\in\mathcal{R}.$ We define
the $\omega-$limit sets $\omega(B)$ and $\omega(\varphi)$ as follows:
\begin{align*}
\omega(B)=  &  \{y\in X:\text{ there are sequences }t_{n}\rightarrow
\infty,y_{n}\in G(t_{n},B)\text{ such that }y_{n}\rightarrow y\},\\
\omega(\varphi)=  &  \{y\in X:\text{ there is a sequence }t_{n}\rightarrow
\infty\text{ such that }\varphi(t_{n})\rightarrow y\}.
\end{align*}
If $\gamma$ is a complete trajectory of $\mathcal{R}$, then the $\alpha-$limit
set is defined by
\[
\alpha(\gamma)=\{y\in X:\text{ there is a sequence }t_{n}\rightarrow
-\infty\text{ such that }\alpha(t_{n})\rightarrow y\}.
\]

Some useful properties of these sets \cite[Lemma 3.4]{ball} are summarized in
the following lemma.

\begin{lemma}
\label{lemadeomegalfa} Assume that $(K1),(K2)$ and $(K4)$ hold. Let $G$ be
asymptotically compact, that is, every sequence $y_{n}\in G(t_{n},B)$, where
$t_{n}\rightarrow\infty$ and $B\subset X$ is bounded, is relatively compact. Then:

\begin{enumerate}
\item For any non-empty bounded set $B, \omega(B)$ is non-empty, compact,
weakly invariant and 
\[
dist_{X}(G(t,B),\omega(B))\rightarrow0, \textit{ as } t\rightarrow+ \infty.
\]

\item For any $\varphi\in\mathcal{R},$ $\omega(\varphi)$ is non-empty,
compact, weakly invariant and 
\[
dist_{X}(\varphi(t),\omega(\varphi
))\rightarrow0, \textit{ as } t\rightarrow+\infty.
\]

\item For any $\gamma\in\mathbb{K},$ $\alpha(\gamma)$ is non-empty, compact,
weakly invariant and 
\[
dist_{X}(\gamma(t),\alpha(\gamma))\rightarrow0, \textit{ as } t\rightarrow-\infty.
\]
\end{enumerate}
\end{lemma}

In order to give a more detailed description of the internal structure of the
attractor under special cases, additional concepts are required.

\begin{definition}
\label{definiciondynamic} Consider a metric space $(X,d)$ and an m-semiflow
$G$.
\end{definition}

\begin{enumerate}
\item We say that $\mathcal{S}=\{\Xi_{1},\ldots,\Xi_{n}\}$ is a family of
isolated weakly invariant sets if there exists $\delta>0$ such that
$\mathcal{O}_{\delta}(\Xi_{i})\cap\mathcal{O}_{\delta}(\Xi_{j})=\emptyset$ for
$1\leq i<j\leq n,$ and each $\Xi_{i}$ is the maximal weakly invariant subset
in $\mathcal{O}_{\delta}(\Xi_{i}):=\{x\in X:dist_{X}(x,\Xi_{i})<\delta\}.$

\item For an m-semiflow $G$ on $(X,d)$ with a global attractor $\mathcal{A}$
and a finite number of weakly invariant sets $\mathcal{S}$, a homoclinic orbit
in $\mathcal{A}$ is a collection $\{\Xi_{p(1)},\ldots,\Xi_{p(k)}\}$
$\subset\mathcal{S}$ and a collection of complete trajectories $\{\gamma
_{i}\}_{1\leq i\leq k}$ of $\mathcal{R}$ in $\mathcal{A}$ such that (putting
$p(k+1):=p(1)$)
\[
\displaystyle\lim_{t\rightarrow-\infty}{dist}_{X}{(\gamma_{i}(t),\Xi_{p(i)}%
)}=0,\displaystyle\lim_{t\rightarrow\infty}{dist}_{X}{(\gamma_{i}%
(t),\Xi_{p(i+1)})}=0,\ 1\leq i\leq k,
\]
and
\[
\text{ for each }i\text{ there exists }t_{i}\in\mathbb{R}\text{ such that
}\gamma_{i}(t_{i})\notin\Xi_{p(i)}\cup\Xi_{p(i+1)}.
\]

\item We say that an m-semiflow $G$ on $(X,d)$ with the global attractor
$\mathcal{A}$ is dynamically gradient if the following two properties
hold:\newline(G1) there exists a finite family $\mathcal{S}=\{\Xi_{1}%
,\ldots,\Xi_{n}\}$ of isolated weakly invariant sets in $\mathcal{A}$ with the
property that any complete trajectory $\gamma$ of $\mathcal{R}$ in
$\mathcal{A}$ satisfies
\[
\displaystyle\lim_{t\rightarrow-\infty}{dist}_{X}{(\gamma(t),\Xi_{i}%
)}=0,\ \lim_{t\rightarrow\infty}{dist}_{X}{(\gamma(t),\Xi_{j})}=0
\]
for some $1\leq i,j\leq n;$\newline(G2) $\mathcal{S}$ does not contain
homoclinic orbits.
\end{enumerate}

\begin{remark}
\label{remark8} In the single-valued case, dynamically gradient semigroups
have been called also gradient-like semigroups \cite{carvalholanga}. Observe
that the above definitions are concerned with weakly invariant families, which
need not to be unitary sets. This is to deal with the more general concept of
generalized gradient-like semigroups \cite{carvalholanga}, in contrast with
gradient-like semigroups (when the invariant sets are unitary).
\end{remark}

Now, we introduce the concept of unstable manifold, that will allow us to
describe more precisely the structure of a global attractor of a dynamically
gradient m-semiflow.

\begin{definition}
Let $G$ be an m-semiflow on a metric space $(X,d)$. The unstable manifold of a
set $\Xi$ is
\[%
\begin{array}
[c]{c}%
W^{u}(\Xi)=\{u_{0} \in X : \text{ there exists complete trajectory }
\gamma\text{ of } \mathcal{R} \text{ such that }\\
\gamma(0)=u_{0} \text{ and } \displaystyle \lim_{t \to-\infty} {dist_X(\gamma
(t),\Xi)}=0 \}.
\end{array}
\]

\end{definition}

Now the following result, relating the global attractor with unstable
manifolds, is standard. The first statement is straightforward to see. The
second one, supposing that the global attractor is compact, follows directly
from the structure described in Theorem \ref{structureattractor} and the
definition of dynamically gradient semiflows.

\begin{lemma}
Consider a metric space $(X,d)$ and a family $\mathcal{R} \subset
\mathcal{C}(\mathbb{R}_{+};X)$ satisfying (K1) and (K2). Suppose that the
associated m-semiflow has a global attractor $\mathcal{A}$. Then, for any
bounded set $\Xi\subset X, W^{u}(\Xi) \subset\bar{\mathcal{A}}$.

Moreover, assume that $\mathcal{R}$ satisfies either (K3) or (K4), and that
the global attractor $\mathcal{A}$ is compact. Suppose also that the
associated m-semiflow $G$ defined in (\ref{defG}) is dynamically gradient.
Then
\[
\mathcal{A}=\bigcup_{i=1}^{n}W^{u}(\Xi_{i}).
\]

\end{lemma}

\section{Robustness of dynamically gradient m-semiflows}

\label{Robustness}

Our first main goal is to prove that a dynamically gradient multivalued
semiflow is stable under perturbations, that is, a family of perturbed
multivalued semiflows remains dynamically gradient if it is close enough to
the original semiflow, generalizing the corresponding result in the
single-valued case \cite{carvalholanga}. This is rigorously formulated in the
following theorem.

\begin{theorem}
\label{teoremaprincipal} Consider a complete metric space $(X,d)$. Let $\eta$
be a parameter in [0,1], $\mathcal{R}_{\eta} \subset\mathcal{C}(\mathbb{R_{+}%
}; X) $ fulfill (K1), (K2), (K3) and (K4), and let $G_{\eta}$ be the
corresponding m-semiflow on X having the global compact attractor
$\mathcal{A}_{\eta}$. Assume that

\begin{itemize}
\item[(H1)] $\overline{\displaystyle \bigcup_{\eta\in[0,1]} \mathcal{A}_{\eta
}}$ \ is compact.

\item[(H2)] $G_{0}$ is a dynamically gradient m-semiflow with finitely many
isolated weakly invariant sets $\mathcal{S}^{0}=\{ \Xi_{1}^{0}, \ldots,
\Xi_{n}^{0} \}$.

\item[(H3)] $\mathcal{A_{\eta}}$ has a finite number of isolated weakly
invariant sets $\mathcal{S_{\eta}}=\{ \Xi_{1}^{\eta}, \ldots, \Xi_{n}^{\eta}
\}$, $\eta\in[0,1]$, which satisfy
\[
\displaystyle \lim_{\eta\to0} {\sup_{1\leq i \leq n} dist_{X}(\Xi_{i}^{\eta},
\Xi_{i}^{0})}=0.
\]

\item[(H4)] Any sequence $\{ \gamma_{\eta} \} $ with $\gamma_{\eta}
\in\mathcal{R_{\eta}}$ such that $\{ \gamma_{\eta}(0) \}$ is converging,
possesses a subsequence $\{ \gamma_{\eta^{\shortmid}} \}$ that converges
uniformly in bounded intervals of $[0, \infty)$ to $\gamma\in\mathcal{R}_{0}$.

\item[(H5)] There exists $\overline{\eta} >0$ and neighborhoods $V_{i}$ of
\ $\Xi_{i}^{0}$ such that \ $\Xi_{i}^{\eta}$ is the maximal weakly invariant
set for $G_{\eta}$ in $V_{i}$ for any $i=1, \ldots, n $ and for each $0<
\eta\leq\overline{\eta}.$
\end{itemize}

Then there exists $\eta_{0}>0$ such that for all $\eta\leq\eta_{0}$, $\{
G_{\eta} \} $ is a dynamically gradient m-semiflow. In particular, the
structure of $\mathcal{A}_{\eta}$ is analogous to that given in
(\ref{attractor}).
\end{theorem}

\begin{proof}
Observe that assumption (H5) concerning certain neighborhood $V_{i}$ of $\Xi_{i}^0$ involves a hyperbolicity condition of $G_{0}$ w.r.t. each $\Xi_{i}^0,$ and as far as (H3) is also assumed, there exist $\{\eta(V_{i})\}_{i=1,\ldots,n}$ such that $\Xi_{i}^\eta\subset V_{i}$ for all $\eta\le\eta(V_{i})$. W.l.o.g. assume that $\delta>0$ is such that $\{x\in X:dist_X(x,\Xi_{i}^0)\le\delta\}\subset V_{i}$ for all $i=1,\ldots, n.$

By Theorem \ref{structureattractor}, we have that $\mathcal A_{\eta}$ is composed by all the orbits of bounded complete trajectories of $\mathcal R_{\eta}$, $\mathbb K_{\eta}.$

We are going to prove by contradiction arguments that there exists $\eta_{0}\in (0,1]$ such that $\{G_{\eta}\}_{\eta\le\eta_{0}}$ is dynamically gradient.

{\bf Step 1:} There exists $\eta_{0}>0$ such that for all $\eta<\eta_{0}$, any bounded complete trajectory $\xi_{\eta}$ of $\mathcal R_{\eta}$ satisfies that there exist $i\in\{1,\ldots ,n\}$ and $t_{0}$ such that for all $t\ge t_{0},$ $dist_X(\xi_{\eta}(t),\Xi_{i}^0)\le\delta$.

After proving the above claim, we consider the sets $B_{\eta}:=\{\xi_{\eta}(s):s\ge t_{0}\} \subset A=\{y :  dist_X(y, \Xi_{i}^{0}) \leq \delta \}$ and $\omega(\xi_{\eta})$. It follows that $\omega(\xi_{\eta})\subset A$, since $dist_X(\xi_{\eta}(t), \omega(\xi_{\eta})) \rightarrow 0$ as $t\rightarrow + \infty$. On the other hand, by Lemma \ref{lemadeomegalfa} $\omega(\xi_{\eta})$ is a weakly invariant set  of $G_{\eta}$ contained in $V_{i}$. By assumption (H5) we have that $\omega({\xi_{\eta}})\subset \Xi_{i}^\eta,$ whence the 'forward part' of property (G1) of a dynamically gradient m-semiflow will follow immediately.

We prove this Step 1 by contradiction. Suppose it does not hold. Then, there exist a sequence $\eta_{k}\to 0$ (as $k\to\infty$) and bounded complete trajectories $\xi_{k}$ of $\mathcal{R}_{\eta_{k}}$ (therefore, from $\mathcal A_{\eta_{k}}$) such that
\begin{equation}
\label{2.1}
\sup_{t\ge t_{0}}dist_X(\xi_{k}(t),\mathcal S^{0})>\delta\ \forall t_{0}\in\mathbb R.
\end{equation}
The set $\{\xi_{k}(0)\}\subset\overline{\cup_{\eta\in[0,1]}\mathcal A_{\eta}}$ is relatively compact from assumption (H1). So, there exists a converging subsequence (relabeled the same) in $X$. From (H4), there exist a subsequence (relabeled the same, again) and $\xi_{0}\in\mathcal R_{0},$ such that $\{\xi_{k}|_{[0,\infty)}\}$ converges to $\xi_{0}$ in bounded intervals of $[0,\infty)$. Actually, if we argue similarly not for time $0,$ but now for times $-1,$ $-2,\ldots,$ and use a diagonal argument, we have that $\xi_{0}=\gamma_{0}|_{[0,\infty)}$ where $\gamma_{0}\in\mathbb K_{0},$ and the convergence of (a subsequence of) $\{\xi_{k}\}$ toward $\gamma_{0}$ holds uniformly in bounded intervals $[a,b]$ of $\mathbb R.$

Since $G_{0}$ is dynamically gradient, there exists $i \in \{1,\ldots ,n\}$ such that $$dist_X(\gamma_{0}(t),\Xi_{i}^0) \rightarrow 0 \textrm{ as } t \rightarrow \infty.$$ Therefore, for all $r\in\mathbb N,$ there exist $t_{r}$ and $k_{r}$ such that $dist_X(\xi_{k}(t_{r}),\Xi_{i}^0)<1/r$ for all $k\ge k_{r}.$ Indeed, this is done as follows: $dist_X(\gamma_{0}(s),\Xi_{i}^0)<1/r$ for all $s\ge t_{r}$ (for some $t_{r},$ w.l.o.g. $t_{r}\ge r>1/\delta$); now, combining this with the uniform convergence on $[0,t_{r}]$ of $\xi_{k}$ toward $\gamma_{0}$, the existence of $k_{r}$ follows.

However, from (\ref{2.1}), there exists $t_{r}'>t_{r}$ such that $dist_X(\xi_{k_{r}}(t),\Xi_{i}^0)<\delta$ for all $t\in[t_{r},t_{r}')$ and $dist_X(\xi_{k_{r}}(t_{r}'),\Xi_{i}^0)=\delta.$

Now we distinguish two cases and we will arrive to the same conclusion in both of them.

{\bf Case (1a):} Suppose that $t_r'-t_r\to\infty$ as $r\to\infty$ (at least for a certain subsequence).

Since $\{\xi_{k_{r}}(t_{r}')\}$ is also relatively compact (by (H1), again), and $\xi_{k_{r}}(t_{r}'+\cdot)$ is a bounded complete trajectory of $\mathcal R_{k_{r}},$ from (H4) we deduce that a subsequence (relabeled the same) is converging on bounded time-intervals of $[0,\infty)$, i.e. $\gamma_{1}(t):=\lim_{r\to\infty}\xi_{k_{r}}(t+t_{r}')$ holds for certain $\gamma_{1}\in\mathcal R_{0}.$ Moreover, as before, a diagonal argument, using not $t_{r}'$ above, but $t_{r}'-1,$ $t_{r}'-2,\ldots$ implies that $\gamma_{1}$ can be extended to the whole real line (the function will still be denoted the same; and the convergence holds in bounded time-intervals of $\mathbb R$), in particular, by (H1) and (H4), $\gamma_{1}\in\mathbb K_{0}$.

Moreover, by its construction, we have that $dist_X(\gamma_{1}(t),\Xi_{i}^0)\le\delta$ for all $t\le 0.$ By Lemma \ref{lemadeomegalfa} we have that the $\alpha$-limit set $\alpha(\gamma_{1})$ is weakly invariant.

As long as $\Xi_{i}^0$ is the biggest weakly invariant set contained in $V_{i},$ we deduce that $dist_X(\gamma_{1}(\tau),\Xi_{i}^0)\rightarrow 0$ when $\tau \rightarrow -\infty$.

On the other hand, from (G1) and (G2) we have that $dist_X(\gamma_{1}(t),\Xi^0_{j})\to 0$ as $t \rightarrow \infty$ for $j\neq i.$

{\bf Case (1b):} Suppose that there exists $C>0$ such that $|t_r'-t_r|\le C$ as $r\to\infty.$ (W.l.o.g. we assume that $t_r'-t_r\to t_*.$)

Recall that $dist_X(\xi_{k_r}(t_r),\Xi^0_i)<1/r.$ By \cite[Lemma 19]{costavalero} $\Xi^0_i$ is closed, so, up to a subsequence $\xi_{k_r}(t_r)\to y\in \Xi^0_i.$ Denote $\xi_{k_r}^1(\cdot )=\xi_{k_r}(\cdot +t_r).$ From (H4), there exist a subsequence $\{\xi_{k_r}^1\}$ and $\xi^1\in\mathcal R_0$ with $\xi^1(0)=y$ such that $\xi_{k_r}^1$ converge towards $\xi^1$ uniformly in bounded intervals of $[0,\infty).$ In particular, $\xi_{k_r}^1(t_r'-t_r)\to\xi^1(t_*),$ so that $dist_X(\xi^1(t_*),\Xi^0_i)\ge\delta.$

Since $\Xi^0_i$ is weakly invariant, there exists $\gamma\in\mathbb K_0$ with $\gamma(0)=\xi^1(0)$ and $\gamma (t)\in\Xi^0_i$ for all $t\in\mathbb R.$ By (K3) consider the concatenation
$$\gamma_1(t):=\left\{\begin{array}{l}\gamma(t),\ \textrm{if }t\le 0,\\
\xi^1(t),\ \textrm{if }t\ge 0.\end{array}\right.$$
Then $\gamma_1\not\equiv\xi^1,$ and by (G1)-(G2) it follows that $dist_X(\gamma_1(t),\Xi^0_j)\to 0$ as $t\to\infty$ with $j\neq i.$ This is exactly the same conclusion we arrived in Case (1a).

Reasoning now with the subsequence $\{\xi_{k_{r}}^1\},$ and proceeding as above, we obtain the existence of $\gamma_{2}\in\mathbb K_{0}$ such that $dist_X(\gamma_{2}(t),\Xi_{j}^0)\to 0$ as $t\to-\infty$ and $dist_X(\gamma_{2}(t),\Xi_{p}^0)\to 0$ as $t\to\infty,$ with $p\not\in\{i,j\}.$

Thus, in a finite number of steps we arrive to a contradiction, since $G_{0}$ satisfies (G2). Therefore, (\ref{2.1}) is absurd, and Step 1 is proved.

{\bf Step 2:} There exists $\eta_{1}>0$ such that for all $\eta<\eta_{1}$, any bounded complete trajectory $\xi_{\eta}$ of $\mathcal R_{\eta}$ satisfies that there exist $j\in\{1,\ldots ,n\}$ and $t_{1}$ such that  $dist_X(\xi_{\eta}(t),\Xi_{j}^0)\le\delta$ for all $t\le t_{1}.$

The above claim can be proved analogously as before, and since for any bounded complete trajectory $\xi_{\eta} \in\mathbb K_{\eta},$ by Lemma \ref{lemadeomegalfa}, $\alpha(\xi_{\eta})$ is weakly invariant for $G_{\eta},$ and contained in some $V_{j},$ the 'backward part' of property (G1) of a dynamically gradient m-semiflow will follow immediately. The same argument is valid for the 'forward part', and so, for all suitable small $\eta,$ $\{G_{\eta}(t):t\ge 0\}$ satisfies (G1).

{\bf Step 3:} There exists $\eta_{2}>0$ such that $\{G_{\eta}\}_{\eta\le\eta_{2}}$ satisfies (G2).

If not, there exist a sequence $\eta_{k}\to 0$, with $G_{\eta_{k}}$ having an homoclinic structure. We may suppose that the number of elements of weakly invariant subsets connected on each homoclinic chain in $\mathcal S_{\eta_{k}}$ is the same. Moreover, by assumption (H3) each $\Xi_{j}^{\eta_{k}}$ is contained in $V_{j}$ for $\eta_{k}$ small enough and w.l.o.g. the order in the route of the homoclinics visiting the $V_{j}$ sets is the same.

Therefore, for $k\ge k_{0}$ there exist a sequence of subsets $\Xi_{p(1)}^{\eta_{k}},\ldots$ $\Xi_{p(l)}^{\eta_{k}}$ in $\mathcal S_{\eta_{k}}$ (with $p(l+1)=p(1)$), and a sequence of complete trajectories $\{\{\xi_{i}^k\}_{i=1}^l\}_{k}$, each collection of $l$ elements in the corresponding attractor $\mathcal A_{\eta_{k}},$ with $$\lim_{t\to-\infty}dist_X(\xi_{i}^k(t),\Xi^{\eta_{k}}_{p(i)})=0,\ \lim_{t\to\infty}dist_X(\xi_{i}^k(t),\Xi^{\eta_{k}}_{p(i+1)})=0, 1\leq i\leq l.$$ If we argue now as in the proof of (G1), we may construct a homoclinic structure of $G_{0},$ getting a contradiction with the fact that the m-semiflow $G_{0}$ is dynamically gradient.
\end{proof}

\begin{remark}
The above result also applies to the particular case of a dynamically gradient
m-semiflow when the weakly invariant families of the original and perturbed
problems are reduced to singletons (Remark \ref{remark8} and \cite[Theorem
1.5]{carvalholanga}).
\end{remark}

\section{An equivalent definition of dynamically gradient families.}

\label{equivalent} We will give an equivalent definition of dynamically
gradient families. For proving the main result in this section we will need a
stronger condition than $(K4)$. Namely, we shall consider the following
stronger condition:

\begin{enumerate}
\item[($\overline{K}4$)] For any sequence $\{ \varphi^{n} \} \subset
\mathcal{R}$ such that $\varphi^{n}(0) \rightarrow\varphi_{0}$ in $X$, there
exists a subsequence $\{\varphi^{n} \}$ and $\varphi\in\mathcal{R}$ such that
$\varphi^{n}$ converges to $\varphi$ uniformly in bounded subsets of $[0,
\infty).$
\end{enumerate}

\begin{remark}
\label{remark2def} We have seen that the property of being dynamical gradient
for a disjoint family of isolated negatively invariant sets $\mathcal{S}%
=\{\Xi_{1}, \ldots, \Xi_{n} \} \subset\mathcal{A}$ is stable under
perturbations. We observe that in the paper \cite{li} a slightly different
definition was used for dynamically gradients families. Namely, instead of
conditions $(G1)$-$(G2)$ it is assumed that any bounded complete trajectory
$\gamma(\cdot)$ satisfies one of the following properties:

\begin{enumerate}
\item $\{\gamma(t) : t \in\mathbb{R} \}\subset\Xi_{i}$ for some $i.$

\item There are $i<j $ for which
\[
\gamma(t) \underset{t\rightarrow\infty}{\rightarrow}\Xi_{i},\ \gamma(t)
\underset{t\rightarrow-\infty}{\rightarrow} \Xi_{j}.
\]

\end{enumerate}
\end{remark}

These assumptions are clearly stronger than $(G1)$-$(G2)$ and imply that the
sets $\Xi_{j}$ are ordered. Our aim is to show that when $\mathcal{S}$ is a
disjoint family of isolated weakly invariant sets, these conditions are
equivalent. For this we will need to introduce the concept of local attractor
and its repeller and study their properties.

We say that $A\subset\mathcal{A}$ is a local attractor in $\mathcal{A}$ if for
some $\varepsilon>0$ we have that $\omega( \mathcal{O}_{\varepsilon}%
(A)\cap\mathcal{A})=A$. Let $A$ be a local attractor in $\mathcal{A}$. Then
its repeller $A^{\ast}$ is defined by
\[
A^{\ast}=\{x\in\mathcal{A}:\omega(x) \backslash A\not = \emptyset\}.
\]

Some properties about local attractors and its repeller as well as the proof
of the following lemmas can be found in \cite{costavalero}.

\begin{lemma}
Assume that $(K1)-(K4)$ hold. Then a local attractor $A$ is invariant.
\end{lemma}

\begin{remark}
Although in \cite{costavalero} the stronger assumption $(\overline{K}4)$ is
assumed, the proof is valid for just $(K4)$.
\end{remark}

\begin{lemma}
\label{lemma5} Assume that $(K1)$-$(K3),\ (\overline{K}4)$ hold and that a
global compact attractor $\mathcal{A}$ exists. Then the repeller $A^{\ast}$ of
a local attractor $A\subset\mathcal{A}$ is weakly invariant and compact.
\end{lemma}

\begin{lemma}
\label{ConvergFunctions}Assume that $(K1)$-$(K3),\ (\overline{K}4)$ hold and
that a global compact attractor $\mathcal{A}$ exists. Let us consider the
sequences $x_{k}\in\mathcal{A},$ $t_{k}\rightarrow+\infty$ and $\varphi_{k}%
($\textperiodcentered$)\in\mathcal{R}$ such that $\varphi_{k}(0)=x_{k}$. Then
from the sequence of maps $\xi_{k}($\textperiodcentered$):[-t_{k}%
,+\infty)\rightarrow\mathcal{A}$ defined by
\[
\xi_{k}(t)=\varphi_{k}(t+t_{k})
\]
one can extract a subsequence converging to some $\psi($\textperiodcentered
$)\in\mathbb{K}$ uniformly on bounded subsets of $\mathbb{R}.$
\end{lemma}

\begin{lemma}
\label{lemma7} Assume that $(K1)$-$(K3),\ (\overline{K}4)$ hold and that a
global compact attractor $\mathcal{A}$ exists. Let $\mathcal{S}=\{\Xi
_{1},\ldots,\Xi_{n}\}\subset\mathcal{A}$ be a disjoint family of isolated
weakly invariant sets. If $G$ is dynamically gradient with respect to
$\mathcal{S}$, then one of the sets $\Xi_{j}$ is a local attractor in
$\mathcal{A}$.
\end{lemma}

\begin{proof}
Let $\delta_{0}>0$ be such that $\mathcal{O}_{\delta_{0}}(\Xi_{i}%
)\cap\mathcal{O}_{\delta_{0}}(\Xi_{j})=\varnothing$ if $i\not =j$ and $\Xi
_{j}$ is the maximal negatively (weakly) invariant set in $\mathcal{O}%
_{\delta_{0}}(\Xi_{j})$ for all $j$. First we will prove the existence of
$j\in\{1,...,n\}$ such that for all $\delta\in(0,\delta_{0})$ there exists
$\delta^{\prime}\in(0,\delta)$ satisfying%
\begin{equation}
\cup_{t\geq0}G(t,\mathcal{O}_{\delta^{\prime}}(\Xi_{j})\cap\mathcal{A}%
)\subset\mathcal{O}_{\delta}(\Xi_{j}). \label{Stab}%
\end{equation}
If not, there would exist $0<\delta<\delta_{0}$ and for each $j$ sequences
$t_{k}^{j}\in\mathbb{R}^{+}$, $x_{k}^{j}\in\mathcal{A}$, $\varphi_{k}^{j}%
\in\mathcal{R}$ with $\varphi_{k}^{j}(0)=x_{k}^{j}$ such that%
\begin{align*}
d(x_{k}^{j},\Xi_{j})  &  <\frac{1}{k},\\
d(\varphi_{k}^{j}(t_{k}^{j}),\Xi_{j})  &  =\delta,\\
d(\varphi_{k}^{j}(t),\Xi_{j})  &  <\delta\text{ for all }t\in\lbrack
0,t_{k}^{j}).
\end{align*}
We have to consider two cases: $t_{k}^{j}\rightarrow+\infty$ or $t_{k}^{j}\leq
C$.

Let $t_{k}^{j}\rightarrow+\infty$. We define the sequence%
\[
\psi_{k}^{j}(t)=\varphi_{k}^{j}(t+t_{k}^{j})\text{ for }t\in\lbrack-t_{k}%
^{j},\infty).
\]
By Lemma \ref{ConvergFunctions} we obtain the existence of a complete
trajectory of $\mathcal{R}$, $\psi^{j}($\textperiodcentered$)$, such that a
subsequence of $\psi_{k}^{j}$ satisfies $\psi_{k}^{j}(t)\rightarrow\psi
^{j}(t)$ for every $t\in\mathbb{R}$. Hence, $d(\psi^{j}(t),\Xi_{j})\leq
\delta<\delta_{0}$ for all $t\leq0$. Therefore, as $\psi^{j}\in\mathbb{K}$,
condition $(G1)$ implies that $d(\psi^{j}(t),\Xi_{j})\rightarrow0$ as
$t\rightarrow-\infty$. On the other hand, since $d(\psi^{j}(0),\Xi_{j}%
)=\delta$, conditions $(G1)-(G2)$ imply that $d(\psi^{j}(t),\Xi_{i}%
)\rightarrow0$ as $t\rightarrow+\infty$, where $i\not =j$.

Let now $t_{k}^{j}\leq C$. We can assume that $t_{k}^{j}\rightarrow t^{j}$. By
$(\overline{K}4)$ we obtain the existence of $\varphi^{j}\in\mathcal{R}$ such
that $\varphi_{k}^{j}$ converges to $\varphi^{j}$ uniformly on bounded sets of
$[0,\infty)$. It is clear then that $d(\varphi^{j}(t^{j}),\Xi_{j})=\delta$. As
$\varphi^{j}(0)\in\Xi_{j}$ and $\Xi_{j}$ is weakly invariant, there exists a
complete trajectory of $\mathcal{R}$, $\psi_{j}^{-}($\textperiodcentered$),$
such that $\psi_{j}^{-}(0)=\varphi^{j}(0)$ and $\psi_{j}^{-}(t)\in\Xi_{j}$ for
all $t\leq0$. Concatenating $\psi_{j}^{-}$ and $\varphi^{j}$ we define%
\[
\psi^{j}(t)=\left\{
\begin{array}
[c]{c}%
\psi_{j}^{-}(t)\text{ if }t\leq0,\\
\varphi^{j}(t)\text{ if }t\geq0,
\end{array}
\right.
\]
which is a complete trajectory by $(K3)$. Again, conditions $(G1)-(G2)$ imply
that $d(\psi^{j}(t),\Xi_{i})\rightarrow0$ as $t\rightarrow+\infty$, where
$i\not =j$.

We have obtained then a connection from $\Xi_{j}$ to a different $\Xi_{i}$.
Since this is true for any $\Xi_{j},$ we would obtain a homoclinic structure,
which contradicts $(G2)$. Therefore, (\ref{Stab}) holds for some $j$. It
follows that
\[
\omega(\mathcal{O}_{\delta^{\prime}}(\Xi_{j})\cap\mathcal{A})\subset
\overline{\mathcal{O}_{\delta}(\Xi_{j})}\subset\mathcal{O}_{\delta_{0}}%
(\Xi_{j}).
\]
Since $\omega(\mathcal{O}_{\delta^{\prime}}(\Xi_{j})\cap\mathcal{A})$ is
weakly invariant, we obtain that $\omega(\mathcal{O}_{\delta^{\prime}}(\Xi
_{j})\cap\mathcal{A})\subset\Xi_{j}$. But $\Xi_{j}\subset G(t,\Xi_{j})\subset
G(t,\mathcal{O}_{\delta^{\prime}}(\Xi_{j})\cap\mathcal{A)}$ for any $t\geq0$
implies the converse inclusion, so that $\Xi_{j}=\omega(\mathcal{O}%
_{\delta^{\prime}}(\Xi_{j})\cap\mathcal{A})$. Thus, $\Xi_{j}$ is a local
attractor in $\mathcal{A}$.
\end{proof}

\bigskip

Now we prove the main result of this section which allows us to establish the
equivalent definition of dynamically gradient families.

\begin{theorem}
\label{resultadogeneral} Assume that (K1)-(K3), $(\overline{K}$4) hold and
that a global compact attractor $\mathcal{A}$ exists. Let $\mathcal{S}%
=\{\Xi_{1},\ldots,\Xi_{n}\}\subset\mathcal{A}$ be a disjoint family of
isolated weakly invariant sets. Then $G$ is dynamically gradient with respect
to $\mathcal{S}$ in the sense of definition \ref{definiciondynamic} if and
only if $\mathcal{S}$ can be reordered in such a way that any bounded complete
trajectory $\gamma(\cdot) $ satisfies one of the following properties:

\begin{enumerate}
\item $\{\gamma(t) : t \in\mathbb{R} \}\subset\Xi_{i}$ for some $i.$

\item There are $i<j$ for which
\[
\gamma( t) \underset{t\rightarrow\infty}{\rightarrow}\Xi_{i},\ \gamma( t)
\underset{t\rightarrow-\infty}{\rightarrow} \Xi_{j}.
\]

\end{enumerate}
\end{theorem}

\begin{proof}
It is obvious that conditions 1-2 imply that $G$ is dynamically gradient. We
shall prove the converse.

By Lemma \ref{lemma7} one of the sets $\Xi_{i}$ is a local attractor. After
reordering the sets, we can say that $\Xi_{1}$ is the local attractor. Let
\begin{equation*}
\Xi_{1}^{\ast}=\{x\in\mathcal{A}:\omega(x)  \backslash \Xi
_{1}\not =\varnothing\}
\end{equation*}
be its repeller, which is weakly invariant by Lemma \ref{lemma5}. Since
$\Xi_{j}$ are closed (cf. \cite[Lemma 19]{costavalero}), weakly invariant and disjoint,
we obtain that $\Xi_{j}\subset\Xi_{1}^{\ast}$ for $j\geq2.$

We will consider only the dynamics inside the repeller $\Xi_{1}^{\ast}$, that is, we define the following set:
\begin{equation*}
\mathcal{R}_{1}=\{\varphi\in\mathcal{R}:\varphi(t)  \in\Xi
_{1}^{\ast}\textrm{ }\forall t\geq0\}.
\end{equation*}
Since $\Xi_{1}^{\ast}$ is weakly invariant, $\mathcal{R}_{1}$ satisfies
$(K1)$. Further, let $\varphi_{\tau}(\cdot)=\varphi(\cdot+\tau)  $, where $\varphi\in\mathcal{R}_{1}$ and $\tau\geq0$. Then it
is clear that $\varphi_{\tau}(t) \in\mathcal{R}_{1}$ for all
$t\geq0$, and then $(K2)$ holds. If $\varphi_{1}(\cdot),\ \varphi_{2}(\cdot)  \in\mathcal{R}_{1}$, it follows by $(K3)$ that the concatenation belongs also to $\mathcal{R}_{1}$. Finally, if
$\varphi_{n}(0)  \rightarrow\varphi_{0}$ with $\varphi_{n}(0)  \in\Xi_{1}^{\ast}$ and $\varphi_{n}(\cdot)  \in\mathcal{R}_{1}$, then $\varphi_{0}\in\Xi_{1}^{\ast}$ (as
$\Xi_{1}^{\ast}$ is closed) and by $(\overline{K}4)  $ passing to a subsequence $\varphi_{n}( t_{n})  \rightarrow \varphi (t)$, for $t_{n}\rightarrow t\geq0$, where $\varphi\in\mathcal{R}$.
Again, the closedness of $\Xi_{1}^{\ast}$ implies that $\varphi\in
\mathcal{R}_{1}$. Hence, $(\overline{K}4)$ also holds. We can
define then the multivalued semiflow $G_{1}:\mathbb{R}^{+}\times\Xi_{1}^{\ast
}\rightarrow P (\Xi_{1}^{\ast}):$
\begin{equation*}
G_{1}(t,x)  =\{y\in\Xi_{1}^{\ast} : y=\varphi(t)
\textrm{ for some } \varphi \in \mathcal{R}_{1} \textrm{, } \varphi(0)
=x\},
\end{equation*}
which is strict by $(K3)  $. This definition is equivalent to the
following one:
\begin{equation*}
\overline{G}_{1}(t,x)=G(t,x)\cap\Xi_{1}^{\ast}\textrm{ for } x \in \Xi_{1}^{\ast}.
\end{equation*}
Indeed, $G_{1}(t,x)  \subset \overline{G}_{1}(t,x)$ is obvious.
Conversely, let $y\in\overline{G}_{1}(t,x)$. Then, $y=\varphi(t)
,\ \varphi(\cdot)  \in\mathcal{R}$, and
$y\in\Xi_{1}^{\ast}$. We state that $\varphi(s)  \in \Xi_{1}
^{\ast}$ for all $0\leq s\leq t$. Assume by contradiction that $\varphi(
s)  \not \in \Xi_{1}^{\ast}$ for $0<s<t$. Therefore, $\omega(
\varphi(s))  \subset\Xi_{1}$. But then by $(K3),$
\begin{equation*}
G(T,y)  \subset G(T,G(t-s,\varphi(s)  ))\subset
G(T+t-s,\varphi(s)  )\rightarrow\Xi_{1}\textrm{ as } T\rightarrow
\infty,
\end{equation*}
which is a contradiction with $y\in\Xi_{1}^{\ast}$. Using again $(
K3)  $ one can define a function $\psi(\cdot)  \in\mathcal{R}_{1}$ such that $\psi(0)=y$, so that $y\in G_{1}(t,x).$

It is clear that $G_{1}$ possesses a global compact attractor, which is the
union of all bounded complete trajectories of $\mathcal{R}_{1}$, and that
$G_{1}$ is dynamically gradient with respect to $\{\Xi_{2},\ldots,\Xi_{n}\}$.
Then, again by Lemma \ref{lemma7} we can reorder the sets in such a way
that $\Xi_{2}$ is a local attractor in $\Xi_{1}^{\ast}$. Let $\Xi_{2,1}^{\ast
}$ be the repeller of $\Xi_{2}$ in $\Xi_{1}^{\ast}$. Then we restrict as
before the dynamics to the set $\Xi_{2,1}^{\ast}$ and so on. Hence, we have
reordered the sets $\Xi_{j}$ in such a way that $\Xi_{1}$ is a local
attractor and $\Xi_{j}$ is a local attractor for the dynamics restricted to
the repeller of the previous local attractor $\Xi_{j-1,j-2}^{\ast}$ for
$j\geq2$, and $\Xi_{i}\subset\Xi_{j-1,j-2}^{\ast}$ if $i\geq j$, where
$\Xi_{1,0}^{\ast}=\Xi_{1}^{\ast}.$

Now, if $\gamma(\cdot)$ is a bounded complete trajectory such that
\begin{equation*}
\gamma(t)  \underset{t\rightarrow\infty}{\rightarrow}\Xi
_{i},\ \gamma(t)  \underset{t\rightarrow-\infty}{\rightarrow}
\Xi_{j},
\end{equation*}
then we shall prove that $i\leq j$. Moreover, if $\gamma(\cdot)$ is not completely contained in some $\Xi_{k}$, then $i<j.$

If $i=1$, then it is clear that $j\geq1$. Also, if there exists $\gamma(
t_{0})  \not \in \Xi_{1}$, then $j>1$, as $\Xi_{1}$ is a local attractor.

Let $i=2$. Then $\gamma(t)  \in\Xi_{1}^{\ast}$ for all
$t\in\mathbb{R}$, and then $\gamma(t)  \underset{t\rightarrow
-\infty}{\rightarrow}\Xi_{1}$ is forbidden. Hence, $j\geq2$. Again, if there
exists $\gamma(t_{0})  \not \in \Xi_{2}$, then the fact that
$\Xi_{2}$ is a local attractor in $\Xi_{1}^{\ast}$ implies that $j>2$.

Further, note that if $i\geq3$, then $\gamma(t) \in \Xi_{1}^{\ast}$ for all $t\in\mathbb{R}$. Also, by induction, it follows that
$\gamma(t) \in \Xi_{k,k-1}^{\ast}$ for all $t \in \mathbb{R}$ and
$2\leq k\leq i-1$. Indeed, let $\gamma(t)  \in\Xi_{k-1,k-2}^{\ast}$ for all $t\in\mathbb{R}$ with $2\leq k\leq i-1$. Then $\gamma(t)  \underset{t\rightarrow\infty}{\rightarrow}\Xi_{i}$ implies clearly
that $\gamma(t) \in \Xi_{k,k-1}^{\ast}$ for all $t\in\mathbb{R}$.
In particular, $\gamma(t) \in \Xi_{i-1,i-2}^{\ast}$ for all
$t\in\mathbb{R}$. Hence, $\Xi_{j}\in\Xi_{i-1,i-2}^{\ast}$, so that $j\geq i$.
Finally, if there exists $\gamma(t_{0})  \not \in \Xi_{i}$, then
$j>i$ as $\Xi_{i}$ is a local attractor in $\Xi_{i-1,i-2}^{\ast}.$
\end{proof}

To finish this section, we recall that the disjoint family of isolated weakly invariant sets
$\mathcal{S}=\{\Xi_{1},\ldots,\Xi_{n}\}\subset\mathcal{A}$ is a Morse
decomposition of the global compact attractor $\mathcal{A}$ if there is a
sequence of local attractors $\emptyset=A_{0}\subset A_{1}\subset\ldots\subset
A_{n}=\mathcal{A}$ such that for every $k\in\{1,\ldots,n\}$ it holds
\[
M_{k}=A_{k}\cap A_{k-1}^{\ast}.
\]

It is well known \cite{li} that for General Dynamical Systems conditions 1-2
in Theorem \ref{resultadogeneral} are equivalent to the fact that
$\mathcal{S}$ generates a Morse decomposition. This fact can be proved also
under conditions (K1)-(K3), $(\overline{K}4)$ \cite{costavalero}.

Thus, Theorem \ref{resultadogeneral} implies that under conditions
$(K1)$-$(K3)$,$(\overline{K}4)$ the family $\mathcal{S}$ generates a Morse
decomposition if and only if $G$ is dynamically gradient.

\section{Application to a reaction-diffusion equation}

\label{application}

We will consider the Chafee-Infante problem%

\begin{equation}
\label{problemainicial}\left\{
\begin{array}
[c]{c}%
\displaystyle \frac{\partial u}{\partial t}-\frac{\partial^{2}u}{\partial
x^{2} }=f(u),\ t>0,\ x\in\left(  0,1\right)  ,\\
u(t,0)=0,\ u(t,1)=0,\\
u(0,x)=u_{0}(x),
\end{array}
\right.
\end{equation}
where $f$ satisfies

\begin{enumerate}
\item[(A1)] $f\in C(\mathbb{R})$;

\item[(A2)] $f(0)=0$;

\item[(A3)] $f^{\prime}\left(  0\right)  >0$ exists and is finite;

\item[(A4)] $f$ is strictly concave if $u>0$ and strictly convex if $u<0;$

\item[(A5)] Growth condition:
\[
|f(u)| \leq C_{1}+C_{2} |u|^{p-1},
\]
where $p \geq2, C_{1},C_{2} >0$;

\item[(A6)] Dissipation condition:

\begin{enumerate}
\item If $p>2$:
\[
f(u)u\leq C_{3}-C_{4}|u|^{p}, \ \ C_{3},C_{4}>0.
\]

\item If $p=2$:
\[
\limsup_{u \rightarrow\pm\infty} \displaystyle \frac{f(u)}{u} \leq0.
\]

\end{enumerate}
\end{enumerate}

\begin{remark}
Note that as a consequence of condition (A6)(b), we have that $f(u)u
\leq(\lambda_{1}-C_{5})u^{2}+C_{6},$ where $C_{5}, C_{6}>0$ and $\lambda
_{1}=\pi^{2}$ is the first eigenvalue of the operator $-\frac{\partial^{2}%
u}{\partial x^{2}}$ with Dirichlet boundary conditions.
\end{remark}

Let $\Omega=(0,1)$ and $1/p+ 1/q=1$. Denote by $(\cdot, \cdot)$ and
$\|\cdot\|_{L^{2}}$ the scalar product and norm in $L^{2}(\Omega)$, by
$\|\cdot\|_{H^{1}_{0}}$ the norm in $H^{1}_{0}(\Omega)$ associated to the
scalar product of gradients in $L^{2}(\Omega)$ thanks to Poincar\'e's
inequality. As usual, let $H^{-1}(\Omega)$ be the dual space to $H_{0}%
^{1}(\Omega)$. Denote by $\langle\cdot, \cdot\rangle$ pairing between the
space $L^{p}(\Omega) \cap H_{0}^{1}(\Omega)$ and its dual $L^{q}(\Omega) \cap
H^{-1}(\Omega)$.

\begin{definition}
The function $u(\cdot) \in C([0,T],L^{2}(\Omega))$ is called a strong solution
of (\ref{problemainicial}) on $[0,T]$ if:

\begin{enumerate}
\item $u(0)=u_{0};$

\item $u(\cdot)$ is absolutely continuous on compact subsets of $(0,T)$;

\item $u(t) \in H^{2}(\Omega) \cap H_{0}^{1}(\Omega)$, $f(u(t)) \in
L^{2}(\Omega)$ for a.e. $t \in(0,T)$ and
\[
\displaystyle \frac{du(t)}{dt}-\Delta u =f(u(t)), \text{ a.e. } t \in(0,T);
\]
where the equality is understood in the sense of the space $L^{2}(\Omega)$.
\end{enumerate}
\end{definition}

\begin{definition}
\label{weakdefinition} The function $u(\cdot) \in C([0,T],L^{2}(\Omega))$ is
called a weak solution of (\ref{problemainicial}) on $[0,T]$ if:

\begin{enumerate}
\item $u \in L^{\infty}(0,T;L^{2}(\Omega))$;

\item $u \in L^{2}(0,T;H^{1}_{0}(\Omega))\cap L^{p}(0,T; L^{p}(\Omega))$;

\item The equality in (\ref{problemainicial}) is understood in the weak sense,
i.e.
\[
\displaystyle \frac{d}{dt} \langle u(t),v \rangle-\langle\Delta u,v
\rangle=\langle f(u(t)),v\rangle,\ \forall v \in H^{1}_{0}(\Omega)\cap
L^{p}(\Omega),
\]
where the equality is understood in the sense of distributions.
\end{enumerate}
\end{definition}

Let us make some comments on the natural relation among the above two
definitions. Let $u(\cdot)$ be a strong solution such that $f(u(\cdot)) \in
L^{2}(0,T;L^{2}(\Omega))$. In view of \cite[Proposition 2.2]{arrietavalero} we
have that $u \in L^{2}(0,T; H^{1}_{0}(\Omega))$, so $\Delta u \in L^{2}(0,T;
H^{-1}(\Omega))$ and then $\frac{du}{dt} \in L^{2}(0,T; H^{-1}(\Omega))$.
Hence, by \cite[Lemma 7.4]{robinson} we get
\[
\displaystyle \langle\frac{du}{dt},v\rangle-\langle\Delta u,v\rangle=\langle
f(u(t)),v \rangle, \ \forall v \in H^{1}_{0}(\Omega).
\]
Using \cite[p.250]{temmam} we obtain
\[
\displaystyle \frac{d}{dt}\langle u,v\rangle-\langle\Delta u,v\rangle=\langle
f(u(t)),v\rangle, \ \forall v \in H^{1}_{0}(\Omega),
\]
so point 3 of Definition \ref{weakdefinition} is satisfied.

Finally, if $p>2$ by condition (A6)(a) we have
\[
|u(t,x)|^{p}\leq\displaystyle\frac{C_{3}}{C_{4}}-\frac{f(u(t,x))u(t,x)}{C_{4}}%
\]
Thus, $f(u)u\in L^{1}((0,T)\times\Omega)$ implies that $u\in L^{p}%
((0,T)\times\Omega)=L^{p}(0,T;L^{p}(\Omega))$. Hence, $u(\cdot)$ is a weak
solution as well.

In view of \cite[p.283]{chepvishik}, for any $u_{0}\in L^{2}(\Omega)$ there
exists at least one weak solution. Moreover, if $f(u(\cdot))\in L^{2}%
(0,T;L^{2}(\Omega))$, then putting $g(\cdot)=f(u(\cdot))$ we obtain by
\cite[p.189]{barbu} that the problem
\[
\left\{
\begin{array}
[c]{lcc}%
\displaystyle\frac{dv}{dt}-\Delta v=g(t) &  & \\
v(0)=u_{0}, &  & \\
&  &
\end{array}
\right.
\]
possesses a unique strong solution $v(\cdot)$. Since this problem has also a
unique weak solution $\tilde{v}(\cdot)$ and the strong solution is a weak
solution as well, then $v(\cdot)=\tilde{v}(\cdot)=u(\cdot)$. Hence $u(\cdot)$
is also a strong solution of problem (\ref{problemainicial}).

Therefore, we have checked that the sets of weak and strong solutions
satisfying $f(u(\cdot))\in L^{2}(0,T;L^{2}(\Omega))$ coincide.

\subsection{Stationary points}

We now focus on the properties of the stationary points. To this end, we have
followed the classic procedure from \cite{henry1981} and \cite{henry1985}.
Moreover, we have also taken some ideas from \cite{simone}.

Let $\mathcal{R} \subset C([0, \infty), L^{2}(\Omega))$ be the set of all weak
solutions of problem (\ref{problemainicial}). Properties $(K1)-(K4)$ are
satisfied [cf. \cite{kapustyankasyanov}], so that a multivalued semiflow is
defined (see Section \ref{preliminaries}). It is shown in \cite[Lemma
12]{kapustyankasyanov} that $v$ is a fixed point of $\mathcal{R}$
(equivalently, of $G$) if and only if $v \in H_{0}^{1}(\Omega)$ and%

\begin{equation}
\label{puntoscriticos}\displaystyle \frac{\partial^{2}v}{\partial x^{2}%
}+f(v)=0, \text{ in } H^{-1}(\Omega).
\end{equation}

The inclusion $H_{0}^{1}(\Omega) \subset L^{\infty}(\Omega)$ implies that
$f(v) \in L^{\infty}(\Omega)$, so that $v \in H^{2}(\Omega) \cap H_{0}%
^{1}(\Omega)$. Therefore, $v(\cdot)$ is a strong solution as well.

Let consider the function $F: \mathbb{R }\rightarrow\mathbb{R}$ defined by
\[
F(s)=\int_{0}^{s} f(r)\mathrm{d}r, \ s \in\mathbb{R}.
\]
We define
\[
a_{-}=inf \{s<0 : \text{ sgn } f(x)= \text{ sgn }x, \ \forall x; s<x<0 \}
\]
and
\[
a_{+}=sup \{s>0 : \text{ sgn } f(x) =\text{ sgn }x, \ \forall x; 0<x<s \}.
\]

If follows from conditions (A2) and (A3) of $f$ that $-\infty\leq
a_{-}<0<a_{+}\leq+\infty$. Since $f$ is positive on $(0,a_{+})$ and negative
on $(a_{-},0)$, we have that $F$ is strictly increasing on $[0,a_{+})$,
strictly decreasing on $(a_{-},0]$ and $F(0)=0$. We consider $E_{+},E_{-}%
\in\lbrack0,\infty]$ defined by
\[
E_{+}=\displaystyle\lim_{s\rightarrow a_{+}}F(s),
\]%
\[
E_{-}=\displaystyle\lim_{s\rightarrow a_{-}}F(s).
\]
Then, $F$ has the inverse functions $U_{+}:[0,E_{+})\rightarrow\lbrack
0,a_{+})$, $U_{-}:[0,E_{-})\rightarrow(a_{-},0]$.

We also define the following functions with domains $(0,E_{+})$ and
$(0,E_{-})$, respectively, with values on $[0,\infty)$:
\[
\tau_{+}(E)=\int_{0}^{U_{+}(E)} \! \!(E-F(u))^{{-1}/{2}}\ \mathrm{d}%
u,\ 0<E<E_{+},
\]
\[
\tau_{-}(E)=\int_{U_{-}(E)}^{0} (E-F(u))^{{-1}/{2}}\ \mathrm{d}u,
\ 0<E<E_{-}.
\]

Let us consider $v_{0}\in\mathbb{R}$ and a solution $u$ of
\begin{align}
\frac{\partial^{2}u}{\partial x^{2}}+f(u)  &  =0,\label{equationfixedpoints}\\
u(0)  &  =0,u^{\prime}(0)=v_{0}.\nonumber
\end{align}
Note that the solution of the problem (\ref{equationfixedpoints}) is unique,
since $f$ is convex for $u<0$ and concave for $u>0$, so it is Lipschitz on
compact intervals \cite[p.4]{daleroberts}, \cite[p.8]{gruber}.

If we define $E=v_{0}^{2}/2$, then:
\[
\frac{(u^{\prime}(x))^{2}}{2}+F(u(x))=E.
\]

On the other hand, the functions $\tau_{+}$, $\tau_{-}$ evaluated in
$E=v_{0}^{2}/2$ give us $\sqrt{2}$ the x-time necessary to go from the initial
condition $u(0)=0$, with initial velocity $v_{0},-v_{0}$ respectively, to the
point where $u^{\prime}(T_{+}(E))=0$. Indeed, $u(x)$ satisfies $\frac
{(u^{\prime}(x))^{2}}{2}+F(u(x))=E$, so $\frac{dx}{du}=\frac{1}{\sqrt{2}}
\frac{1}{\sqrt{E-F(u)}}$. Since $u^{\prime}(T_{+}(E))=0$ for $u=U_{+}(E)$,
then
\[
\sqrt{2} \int_{0}^{T_{+}(E)}\!\!\!\!\!\!\! 1 \ \ \mathrm{d}x=\int_{0}%
^{U_{+}(E)}\frac{1}{\sqrt{E-F(u)}} \mathrm{d}u=\tau_{+}(E).
\]

By symmetry with respect to the $u$ axis, the $x-$time it takes for $u(x)$ to
go from $(U^{+}(E),0)$ to $(0,-v_{0})$ is $T_{+}(E)$. By this way, if
$2T_{+}(E)=1$, that is, $\tau^{+}(E)=\frac{1}{\sqrt{2}}$, then $u(\cdot)$ is a
solution satisfying the boundary conditions $u(0)=u(1)=0$. Applying a similar
reasoning for $\tau^{-}(E)$, we obtain that $u$ satisfies the boundary
conditions if, and only if, $E$ satisfies for some $k \in\mathbb{N} $ only one
of the following conditions:
\begin{equation}
\label{primera}k\tau_{+}(E)+(k-1)\tau_{-}(E)=\frac{1}{\sqrt{2}},
\end{equation}

\begin{equation}
\label{segunda}k\tau_{-}(E)+(k-1)\tau_{+}(E)=\frac{1}{\sqrt{2}},
\end{equation}

\begin{equation}
\label{tercera}k\tau_{+}(E)+k\tau_{-}(E)=\frac{1}{\sqrt{2}}.
\end{equation}

\begin{remark}
\label{numzeros} Note that if $E$ satisfies (\ref{primera}) or (\ref{segunda})
for a certain $k$, then $u$ has $2k$ zeros and if $E$ satisfies (\ref{tercera}%
), then $u$ has $2k+1$ zeros. Our goal is to solve these equations for $E$ as
a function of $f^{\prime}(0)$. To this end, we study the properties of
$\tau_{\pm}$.
\end{remark}

In order to obtain solutions of the equations (\ref{primera}), (\ref{segunda})
and (\ref{tercera}) it is necessary to make a change of variable for the
functions $\tau_{\pm}$. Given $E \in(0,E_{\pm})$, we put
\[
Ey^{2}=F(u), \ 0 \leq y \leq1, \ 0 \leq u \leq U_{+}(E)
\]
and
\[
Ey^{2}=F(u), \ -1 \leq y \leq0, \ U_{-}(E) \leq u \leq0.
\]
Hence, $du=(2yE/f(u))dy$ and $E-F(u)=E(1-y^{2})$. By this change, we obtain
\[
\tau_{+}(E)=2 \sqrt{E} \displaystyle \int_{0}^{1} (1-y^{2})^{-1/2}\frac
{y}{f(u)}dy, 0<E<E_{+}; \ u=U_{+}(Ey^{2}), 0 \leq y \leq1;
\]
\[
\tau_{-}(E)=2 \sqrt{E} \displaystyle \! \int_{-1}^{0} (1-y^{2})^{-1/2}\frac
{y}{f(u)}dy, 0<E<E_{-}; u=U_{-}(Ey^{2}), -1 \leq y \leq0.
\]

The next results show some properties of these functions.

\begin{theorem}
The functions $\tau_{\pm}$ satisfy
\[
\displaystyle\lim_{E \rightarrow0^{+}} \tau_{\pm}(E)=\frac{\pi}{(2f^{\prime
}(0))^{{1}/{2}}}.
\]

\end{theorem}

\begin{proof}
Since $f^{\prime}(0)>0$ and $f(0)=0$, given $\varepsilon\in(0,1)$, there
exists $\delta>0$ such that
\begin{equation}%
\begin{array}
[c]{c}%
f^{\prime}(0)(1-\varepsilon)u\leq f(u)\leq f^{\prime}(0)(1+\varepsilon
)u,\ \ 0\leq u\leq\delta.\\
\displaystyle\frac{1}{f^{\prime}(0)(1+\varepsilon)}\leq\frac{u}{f(u)}\leq
\frac{1}{f^{\prime}(0)(1-\varepsilon)},\ 0\leq u\leq\delta.
\end{array}
\label{continu1}%
\end{equation}
Moreover, as $U_{+}(E)$ is continuous at $0$, given $\delta>0$, there exists
$\eta>0$ such that for $0<E\leq\eta,\ U_{+}(E)\leq\delta$. Now, if we
integrate (\ref{continu1}) between $0$ and $u$ we obtain the following
inequality
\[
\displaystyle\frac{f^{\prime}(0)}{2}(1-\varepsilon)u^{2}\leq F(u)\leq
\frac{f^{\prime}(0)}{2}(1+\varepsilon)u^{2},\ 0\leq u\leq\delta.
\]
Using the change of variable $Ey^{2}=F(u)$, we have
\[%
\begin{array}
[c]{c}%
\left(  \displaystyle\frac{f^{\prime}(0)(1-\varepsilon)}{2E}\right)  ^{{1}%
/{2}}u\leq y\leq\left(  \displaystyle\frac{f^{\prime}(0)(1+\varepsilon)}%
{2E}\right)  ^{{1}/{2}}u,\text{ \ for \ }0<E\leq\eta,\ 0\leq y\leq1.
\end{array}
\]
Dividing the previous expression by $f(u)$ and using (\ref{continu1}) we
obtain
\[%
\begin{array}
[c]{c}%
\left(  \displaystyle\frac{1-\varepsilon}{2Ef^{\prime}(0)(1+\varepsilon)^{2}%
}\right)  ^{{1}/{2}}\!\!\!\leq\displaystyle\frac{y}{f(u)}\leq\left(
\displaystyle\frac{1+\varepsilon}{2Ef^{\prime}(0)(1-\varepsilon)^{2}}\right)
^{{1}/{2}}\!\!\!,\!\!\!\text{ \ for \ }0<E\leq\eta,\ 0\leq y\leq1.
\end{array}
\]
Now if we multiply by $2\sqrt{E}(1-y^{2})^{-\frac{1}{2}}$ and integrate from
$0$ to $1$, we get
\[%
\begin{array}
[c]{c}%
\displaystyle\!\!\!\pi\left(  \frac{1-\varepsilon}{2f^{\prime}%
(0)(1+\varepsilon)^{2}}\right)  ^{{1}/{2}}\leq\tau_{+}(E)\leq\pi\left(
\frac{1+\varepsilon}{2f^{\prime}(0)(1-\varepsilon)^{2}}\right)  ^{{1}/{2}%
},\!\!\!0<E\leq\eta.
\end{array}
\]
Finally, taking $\varepsilon\rightarrow0$, the theorem follows. The proof for
$\tau_{-}$ is analogous.
\end{proof}

\begin{theorem}
The functions $\tau_{\pm}$ are strictly increasing on their domains.
\end{theorem}

\begin{proof}
Let consider the expression of $\tau_{+}$ and $0<E_{1}<E_{2}<E_{+}$. Then,
\[
\tau_{+}(E_{2})-\tau_{+}(E_{1})=\displaystyle\int_{0}^{1}\frac{2y}%
{\sqrt{1-y^{2}}}\left[  \displaystyle\frac{\sqrt{E_{2}}}{f(U^{+}(E_{2}y^{2}%
))}-\frac{\sqrt{E_{1}}}{f(U^{+}(E_{1}y^{2}))}\right]  dy.
\]

From \cite[p.8]{gruber} we have that the function $f$ is differentiable almost
everywhere in $\mathbb{R}$, so $\alpha(E)=\displaystyle\frac{\sqrt{E}}%
{f(U^{+}(Ey^{2}))}$ is differentiable as well. Hence,
\[
\alpha^{\prime}(E)=\displaystyle\frac{f^{2}(U^{+}(Ey^{2}))-2y^{2}Ef^{\prime
}(U^{+}(Ey^{2}))}{2\sqrt{E}f^{3}(U^{+}(Ey^{2}))}.
\]
Recall the change of variable $F(u)=Ey^{2}$. Consider the numerator of
$\alpha^{\prime}$, that is $\beta(u)=f^{2}(u)-2F(u)f^{\prime}(u)$, then we
obtain
\[
\beta(u)=\displaystyle2\int_{0}^{u}f(s)(f^{\prime}(s)-f^{\prime}%
(u))ds,\ 0<s<u.
\]
Since $f$ is strictly concave, if $s<u,$ then $f^{\prime}(s)>f^{\prime}(u)$
(cf. \cite[p.5]{daleroberts}). As a result, $\beta(u)>0$.

In order to finish
the proof rigorously, we have to justify the previous calculations. Indeed,
from \cite[p.5]{gruber}, we have that the function $f$ is absolutely
continuous and from \cite[p.16]{barbu}, $f^{\prime}\in L_{loc}^{1}$.
Therefore, $\alpha^{\prime}\in L_{loc}^{1}$ and $\alpha^{\prime}>0$ a.e.,
which implies that $\alpha(E)$ is strictly increasing and the proof is
finished.

 The claim for $\tau_{-}(E)$ follows analogously.
\end{proof}

\begin{theorem}
The functions $\tau_{\pm}$ satisfy
\[
\displaystyle\lim_{E \rightarrow E^{\pm}} \tau_{\pm}(E)=\infty
\]
Then, $\tau_{\pm}:(0,E^{\pm}) \rightarrow\left(  \displaystyle \frac{\pi
}{(2f^{\prime}(0))^{{1}/{2}}},\infty\right)  $.
\end{theorem}

\begin{proof}
\textbf{Case $a_{+}<\infty$.} Then, we have $f(a_{+})=0$ and $\bar{u}%
(x)=a_{+}$ is a constant solution to the problem $\frac{\partial^{2}%
u}{\partial x^{2}}+f(u)=0$. Let us consider $E_{+}=F(a_{+})$ and the solution
$u$ to this problem satisfying the conditions $u(0)=0,u^{\prime}%
(0)=v_{0},E=\frac{1}{2}v_{0}^{2}.$ As $a_{+}$ is a constant solution, by
uniqueness $\tau_{+}(E^{+})=\infty$. Therefore, given $T>0$, there exists
$\delta>0$ such that if $E>E_{+}-\delta$, then $\tau_{+}(E)>T$, which follows
from the continuity of $u$ with respect to its initial conditions.

\textbf{Case $a_{+}=\infty.$} Note that if $p>2$, then $a_{+}<\infty$.
Therefore, $p=2$. In this case, $f(u)>0$ for all $u\in(0,\infty)$. From
condition (A5), there exist $\alpha,\beta>0$ such that $f(u)\leq\alpha+\beta
u$. For $u>0$ we have
\[
\frac{f(u)}{u^{2}}\leq\frac{\alpha}{u^{2}}+\frac{\beta}{u}.
\]
Hence, $f(u)/u^{2}\rightarrow0,$ as $u\rightarrow\infty$.

On the other hand,
$\int_{0}^{u}f(s)ds\leq\int_{0}^{u}(\alpha+\beta s)\ ds$. Thus, we have
$\displaystyle F(u)\leq\alpha u+\beta{u^{2}}/{2}$ and
\[
\displaystyle0\leq\frac{F(u)}{u^{3}}\leq\frac{\alpha}{u^{2}}+\frac{\beta}%
{2}\frac{1}{u}.
\]
Hence, $F(u)/u^{3}\rightarrow0,$ as $u\rightarrow\infty$.

We claim that \textbf{$\displaystyle\lim_{u\rightarrow0^{+}}f(u)/u^{2}%
=\infty.$} Indeed, since $f^{\prime}(0)$ exists, for any $\varepsilon
\in(0,f^{\prime}(0)),$ there exists $\delta>0$ such that $|f^{\prime
}(0)-f(u)/u|<\varepsilon$, for any $|u|<\delta$. Thus, dividing by $u^{2}$, we
obtain
\[
\displaystyle\frac{u(f^{\prime}(0)-\varepsilon)}{u^{2}}<\frac{f(u)}{u^{2}%
}<\frac{u(\varepsilon+f^{\prime}(0))}{u^{2}}%
\]
and the result follows.

Since $f(u)/u^{2}\rightarrow0$, as $u\rightarrow
\infty$, and $f(u)/u^{2}\rightarrow\infty$, as $u\rightarrow0^{+}$, for any
$\varepsilon>0$, there exists a first value $u_{0}\in(0,\infty)$ where
$f(u_{0})/u_{0}^{2}=\varepsilon.$ Hence,
\[
\frac{f(u)}{u^{2}}>\varepsilon,\ 0<u<u_{0}.
\]

From the above expression, we have $\int_{0}^{u}f(s)ds>\int_{0}^{u}\varepsilon
s^{2}ds$ and $\varepsilon u^{3}/3<F(u)$. Then, $F(u)/u^{3}>\varepsilon/3$, if
$0<u\leq u_{0}$. Since $F(u)/u^{3}\rightarrow0$, as $u\rightarrow\infty$, we
deduce that there exists $\overline{u}>u_{0}$ such that $F(\overline
{u})/\overline{u}^{3}=\varepsilon/3$. Hence, we have
\[
\displaystyle\frac{F(u)}{u^{3}}>\frac{\varepsilon}{3},\ 0<u<\overline{u},
\]
with $F(\overline{u})=\frac{\varepsilon}{3}\overline{u}^{3}$.

Now, computing
$\tau_{+}$ in $\overline{E}=F(\overline{u})$, we have
\begin{equation}%
\begin{array}
[c]{c}%
\tau_{+}(\overline{E})=\displaystyle\int_{0}^{U^{+}(\overline{E}%
)}\displaystyle\frac{1}{\sqrt{\overline{E}-F(u)}}du=\displaystyle\int%
_{0}^{\overline{u}}\displaystyle\frac{1}{\sqrt{\frac{\varepsilon}{3}%
\overline{u}^{3}-F(u)}}du\nonumber\\
\\
\geq\displaystyle\int_{0}^{\overline{u}}\displaystyle\frac{1}{\sqrt
{\frac{\varepsilon}{3}\overline{u}^{3}-\frac{\varepsilon}{3}u^{3}}%
}du=\displaystyle\frac{\sqrt{3}}{\sqrt{\varepsilon}}\displaystyle\int%
_{0}^{\overline{u}}\displaystyle\frac{1}{\sqrt{\overline{u}^{3}-u^{3}}}du\\
\\
=\displaystyle\frac{\sqrt{3}}{\sqrt{\varepsilon}}\displaystyle\int_{0}%
^{1}\displaystyle\frac{\overline{u}}{\sqrt{\overline{u}^{3}-\overline{u}%
^{3}t^{3}}}dt=\displaystyle\frac{\sqrt{3}}{\sqrt{\varepsilon}}%
\displaystyle\frac{\overline{u}}{\sqrt{\overline{u}^{3}}}\int_{0}^{1}\left(
1-t^{3}\right)  ^{-\frac{1}{2}}dt\\
\\
=\displaystyle\frac{\sqrt{3}}{\sqrt{\varepsilon}}\displaystyle\frac
{\overline{u}}{\sqrt{\overline{u}^{3}}}\frac{1}{3}\int_{0}^{1}s^{\frac{1}%
{3}-1}\left(  1-s\right)  ^{\frac{1}{2}-1}ds\\
\\
=\displaystyle\frac{1}{\overline{u}^{\frac{1}{2}}}\frac{1}{\sqrt{\varepsilon}%
}\frac{\sqrt{3}}{3}\mathcal{B}\left(  \frac{1}{2},\frac{1}{3}\right)  .
\end{array}
\end{equation}

Recall that $\varepsilon\overline{u}^{3}=3F(\overline{u})$. Then,
\[
\varepsilon\overline{u}=3\displaystyle\frac{F(\overline{u})\overline{u}%
}{\overline{u}^{3}}=\displaystyle3\frac{F(\overline{u})}{\overline{u}^{2}}%
\]
Taking $\varepsilon\rightarrow0$, by construction $\overline{u}\rightarrow
\infty$. Therefore, from condition (A6)(b) we have that $\lim_{u\rightarrow
\infty}f(u)/u\leq0$, so the last expression tends to $0$ and $\tau
_{+}(\overline{E})\rightarrow\infty$.
\end{proof}

\begin{theorem}
\label{teorema18} Consider
\[
\lambda_{n}=n^{2}\pi^{2}.
\]
Then, for each $n\geq1$, there exist two continuous functions $E_{n}^{\pm}:
[\lambda_{n},\infty) \rightarrow[0,E_{\pm})$ with the following properties:

\begin{enumerate}
\item For each integer $k \geq1$ and for $f^{\prime}(0) \in[\lambda
_{2k-1},\infty)$ the only solution of the equation (\ref{primera}) (resp.
\ref{segunda}) is the value $E^{+}_{2k-1}(f^{\prime}(0))$ (resp. $E^{-}%
_{2k-1}(f^{\prime}(0))$);

\item For each integer $k\geq1$ and for $f^{\prime}(0)\in\lbrack\lambda
_{2k},\infty)$ the only solution of the equation (\ref{tercera}) is the value
$E_{2k}^{-}(f^{\prime}(0))=E_{2k}^{+}(f^{\prime}(0))=E_{2k};$

\item For each integer $n \geq1$, $E^{\pm}_{n}(f^{\prime}(0))=0$, if
$f^{\prime}(0)=\lambda_{n}$.
\end{enumerate}
\end{theorem}

\begin{proof}
Let be $n\geq1$. If $n$ is odd, then $n=2k-1$ for $k\geq1$. First, we prove
that we can define the function
\begin{equation}%
\begin{array}
[c]{c}%
E_{n}^{\pm}:[\lambda_{n},\infty)\longrightarrow\lbrack0,E_{\pm})\nonumber
\end{array}
\end{equation}
by putting $E_{n}^{\pm}(f^{\prime}(0))=E$, where $E$ satisfies $k\tau_{\pm
}(E)+(k-1)\tau_{\mp}(E)=1/\sqrt{2}$.

Consider the function
\[
h_{\pm}^{n}:(0,E_{\pm})\longrightarrow(n\pi/\sqrt{2f^{\prime}(0)},\infty),
\]
defined by $h_{\pm}^{n}(E):=k\tau_{\pm}(E)+(k-1)\tau_{\mp}(E)$. If $f^{\prime
}(0)>\lambda_{n}$ then, as $h_{\pm}$ is a strictly increasing function, there
exists a unique $E_{2k-1}^{\pm}\in(0,E_{\pm})$ such that $h_{\pm}^{n}%
(E_{2k-1}^{\pm})=1/\sqrt{2}$.

Since $h_{\pm}$ has inverse, $E_{2k-1}^{\pm
}=(h_{\pm}^{n})^{-1}({1}/\sqrt{2})$ is the solution of the expressions
(\ref{primera}) and (\ref{segunda}). Moreover, $E_{2k-1}^{\pm}(\lambda_{n})=0$
by construction.

Second, if $n$ is even, then $n=2k$ for $k\geq1$. As before, we consider
$h_{\pm}^{n}(E):=k\tau_{\pm}(E)+k\tau_{\mp}(E)$. Since it is an increasing
function, for $f^{\prime}(0)>\lambda_{n}$, there exists a unique $E_{2k}%
\in(0,E_{\pm})$ such that $h_{\pm}^{n}(E_{2k})=1/\sqrt{2}$. Analogously, we
obtain the solution of the expression (\ref{tercera}), $E_{2k}^{\pm}=(h_{\pm
}^{n})^{-1}({1}/\sqrt{2})$, and $E_{2k-1}^{\pm}(\lambda_{n})=0$.
\end{proof}

\begin{theorem}
\label{teorema19}For each $n\geq0$ and $f^{\prime}(0)\in\lbrack\lambda
_{n},\infty)$, the equation (\ref{puntoscriticos}) has two new more solutions
$v_{n}^{\pm}$ with the following properties:

\begin{enumerate}
\item $a_{-}<u_{n}^{\pm}(x)<a_{+}$ for all $x\in\lbrack0,1];$

\item If $f^{\prime}(0)=\lambda_{n}$, then $v_{n}^{\pm}=0$;

\item For $f^{\prime}(0)\in(\lambda_{n},\infty)$, $v_{n}^{\pm}$ has $n+1$
zeros in $[0,1]$. Denoting these zeros by $x_{q}^{\pm}$, $q=0,1,\ldots,n$ with
$0=x_{0}^{\pm}<x_{1}^{\pm}<x_{2}^{\pm}<\ldots<x_{n}^{\pm}=1,$ we have
$(-1)^{q}v_{n}^{+}(x)>0$ for $x_{q}^{+}<x<x_{q+1}^{+},q=0,1,\ldots,n-1$ and
$(-1)^{q}v_{n}^{-}(x)<0$ for $x_{q}^{-}<x<x_{q+1}^{-},q=0,1,\ldots,n-1$. Also,
$v_{n}^{+}=-v_{n}^{-},$ if $f$ is odd;
\end{enumerate}
\end{theorem}

\begin{proof}
The first point follows from $F(u_{n}^{\pm}(x))\leq E<E_{\pm}.$

The second point follows from the third one of Theorem \ref{teorema18}.
Indeed, for each $n\geq1$ and $f^{\prime}(0)\in\lbrack\lambda_{n},\infty)$ we
have the values $E_{n}^{\pm}(f^{\prime}(0))$ by the above theorem. Also, we
have a solution of the equation (\ref{puntoscriticos}) which is denoted by
$v_{n}^{\pm}$. If $f^{\prime}(0)=\lambda_{n}$, then $E_{n}^{\pm}(\lambda
_{n})=0$ and $v_{0}=0$, so $v_{n}^{\pm}=0$.

The third point follows by Remark \ref{numzeros}. If $f$ is odd, then
$-U^{-}(E)=U^{+}(E)$, $\tau_{+}(E)=\tau_{-}(E)$, so we have $v_{n}^{+}%
=-v_{n}^{-}$.
\end{proof}

\begin{corollary}
If $n^{2}\pi^{2}<f^{\prime}(0)\leq(n+1)^{2}\pi^{2},n\in\mathbb{N}$, then there
are $2n+1$ fixed points: $0,\ v_{1}^{\pm},\ ...,v_{n}^{\pm}$, where
$v_{j}^{\pm}$possesses $j+1$ zeros in $[0,1].$
\end{corollary}

\subsection{Approximations}

From now on, we shall consider the following family of Chafee-Infante
equations
\begin{equation}
\label{aproximaciones}\left\{
\begin{array}
[c]{c}%
\displaystyle \frac{\partial u}{\partial t}-\frac{\partial^{2}u}{\partial
x^{2} }=f_{\varepsilon}(u),\ t>0,\ x\in\left(  0,1\right)  ,\\
u(t,0)=0,\ u(t,1)=0,\\
u(0,x)=u_{0}(x),
\end{array}
\right.
\end{equation}
where $\varepsilon\in(0,1]$ is a small parameter and $f_{\varepsilon}$ satisfies

\begin{itemize}
\item[($\widetilde{A1}$)] $f_{\varepsilon}\in C(\mathbb{R})$ and\ is non-decreasing;

\item[($\widetilde{A2}$)] $f_{\varepsilon}(0)=0$;

\item[($\widetilde{A3}$)] $f_{\varepsilon}^{\prime}\left(  0\right)  >0 $
exists, is finite, monotone in $\varepsilon$ and $f_{\varepsilon}^{\prime
}\left(  0\right)  {\rightarrow}\infty$, as ${\varepsilon\rightarrow0^{+}}$;

\item[($\widetilde{A4}$)] $f_{\varepsilon}$ is strictly concave if $u>0$ and
strictly convex if $u<0;$

\item[($\widetilde{A5}$)] $-1<f_{\varepsilon}\left(  s\right)  <1,$ for all
$s,$ and
\begin{equation}
\label{adecuacion}|f_{\varepsilon}(s) -H_{0}(s)| <\varepsilon, \ \text{ if }
|s|>\varepsilon,
\end{equation}
where%

\[
H_{0}(u)= \left\{
\begin{array}
[c]{lcc}%
-1, & \text{if} & u< 0,\\
\left[  -1,1\right]  , & \text{if} & u=0,\\
1, & \text{if} & u>0,
\end{array}
\right.
\]
is the Heaviside function.
\end{itemize}

Conditions (A1)-(A6) are satisfied with $p=2,$ so problem
(\ref{aproximaciones}) is a particular case of (\ref{problemainicial}).

Our aim now is to prove that for $\varepsilon$ sufficiently small the
multivalued semiflow $G_{\varepsilon}$ generated by the weak solutions of
problem (\ref{aproximaciones}) is dynamically gradient. Problem
(\ref{aproximaciones}) is an approximation of the following problem, governed
by a differential inclusion
\begin{equation}
\left\{
\begin{array}
[c]{c}%
\displaystyle\frac{\partial u}{\partial t}-\frac{\partial^{2}u}{\partial
x^{2}}\in H_{0}(u),\text{ on }\Omega\times(0,T),\\
u|_{\partial\Omega}=0,\\
u(0,x)=u_{0}(x).
\end{array}
\right.  \label{inclusion1}%
\end{equation}

We say that the function $u \in C([0,t], L^{2}(\Omega))$ is a strong solution
of (\ref{inclusion1}) if

\begin{enumerate}
\item $u(0)=u_{0}$;

\item $u(\cdot)$ is absolutely continuous on $(0,T)$ and $u(t) \in
H^{2}(\Omega) \cap H_{0}^{1}(\Omega)$ for a.e. $t \in(0,T)$;

\item There exists a function $g(\cdot)$ such that $g(t) \in L^{2}(\Omega)$,
a.e. on $(0,T)$, $g(t,x) \in H_{0}(u(t,x))$, for a.e. $(t,x) \in(0,T)
\times\Omega$, and
\[
\displaystyle \frac{du}{dt}- \frac{\partial^{2}u}{\partial x^{2}}-g(t)=0,
\text{ a.e. } t \in(0,T).
\]

\end{enumerate}

In this case we put $\mathcal{R}$ as the set of all strong solutions such that
the map $g$ belongs to $L^{2}(0,T;L^{2}(\Omega))$. Conditions (K1)-(K4) are
satisfied (cf. \cite{costavalero}) and the map $G:\mathbb{R}_{+}\times
L^{2}(\Omega)\to P(L^{2}(\Omega))$ defined by (\ref{defG}) is a strict
multivalued semiflow possessing a global compact attractor $\mathcal{A}_{0}$
(cf. \cite{valero2001}) in $L^{2}(\Omega),$ which is connected (cf.
\cite{valero2005}). The structure of this attractor is studied in
\cite{arrietavalero}. It is shown that there exists an infinite (but
countable) number of fixed points
\[
v_{0}=0, v_{1}^{+}, v_{1}^{-}, \ldots, v_{n}^{+}, v_{n}^{-}, \ldots,
\]
and that $\mathcal{A}_{0}$ consists of these fixed points and all bounded
complete trajectories $\psi(\cdot)$, which always connect two fixed points,
that is,
\begin{equation}
\label{trayectoriasconectando}%
\begin{array}
[c]{l}%
\psi(t) \rightarrow z_{1} \text{ as } t \rightarrow\infty,\\
\psi(t) \rightarrow z_{2} \text{ as } t \rightarrow-\infty,
\end{array}
\end{equation}
where $z_{i}=0, z_{i}=v_{n}^{+}$ or $z_{i}=v_{n}^{-}$ for some $n \geq1$.
Moreover, if $\psi$ is not a fixed point, then either $z_{2}=0$ and
$z_{1}=v_{n}^{\pm}$, for some $n \geq1$, or $z_{2}=v_{k}^{\pm}$, $z_{1}%
=v_{n}^{\pm}$ with $k>n$.

Denote
\[
Z_{n}=\left(  \cup_{k \geq n} \{v_{k}^{\pm} \} \right)  \cup\{v_{0} \}
\]
and define the sets
\[
\Xi_{k}^{0}= \{v_{k}^{+}, v_{k}^{-} \}, \ \ 1 \leq k \leq n-1,
\]
\[
\Xi_{n}^{0}= \left\{
\begin{array}
[c]{ccc}%
y : \exists\psi\in\mathbb{K} \text{ such that (\ref{trayectoriasconectando})
holds} \text{ with } z_{j} \in Z_{n}, \  &  & \\
j=1,2 \text{ and } y=\psi(t) \text{ for some } t \in\mathbb{R} &  &
\end{array}
\right\}  ,
\]
where $\mathbb{K}$ stands for the set of all bounded complete trajectories. We
note that set $\Xi_{n}^{0}$ contains the fixed points in $Z_{n}$ and all
bounded complete trajectories connecting them.

\begin{remark}
\label{condicion2delmaintheorem} It is known \cite{costavalero} that the
family $\mathcal{M}=\{\Xi_{1}^{0}, \ldots, \Xi_{n}^{0} \}$ is a disjoint
family of isolated weakly invariant sets and that $G_{0}$ is dynamically
gradient with respect to $\mathcal{M}$ in the sense of Remark \ref{remark2def}%
. Since assumption ($\overline{K}$4) holds true (see \cite[Lemma
31]{costavalero}), Theorem \ref{resultadogeneral} implies that $G_{0}$ is
dynamically gradient with respect to $\mathcal{M}$ in the sense of Definition
\ref{definiciondynamic}.
\end{remark}

Now our purpose is to adapt some lemmas from \cite[p.2979]{arrietavalero} to
problem (\ref{aproximaciones}). In view of Theorems \ref{teorema18} and
\ref{teorema19} and the third condition on $f_{\varepsilon}$, there exists a
sequence $\overline{\varepsilon}_{k} \rightarrow0$, as $k\rightarrow\infty$,
such that for every $\varepsilon\in(\overline{\varepsilon}_{k},\overline
{\varepsilon}_{k+1}]$ and any $k \geq1$ problem (\ref{aproximaciones}) has
exactly $2k+1$ fixed points $\{v_{0}^{\varepsilon}=0, \{v_{\varepsilon,j}%
^{+}\}_{j=1}^{k} \}$ such that for each $1 \leq n \leq k$ $v_{\varepsilon
,n}^{\pm}$ has $n+1$ zeros in $[0,1]$.

Let us consider a sequence $\{\varepsilon_{m}\}$ converging to zero.

\begin{lemma}
\label{lemaconvergenciaA0} Let $n \in\mathbb{N}$ be fixed. Then,
$v_{\varepsilon_{m},n}^{+}$ (resp. $v_{\varepsilon_{m},n}^{-}$) do not
converge to $0$ in $H^{1}_{0}(0,1)$ as $\varepsilon_{m} \rightarrow0$.
\end{lemma}

\begin{proof}
Suppose that $v_{\varepsilon_{m},n}^{+}\rightarrow0$ in $H_{0}^{1}(0,1)$. Then
$v_{\varepsilon_{m},n}^{+}\rightarrow0$ in $C([0,1]).$ By Remark
\ref{numzeros}, $v_{\varepsilon_{m},n}^{+}$ has a unique maximum in
$a\in(0,x_{1}^{+})$ and by the properties of $\tau_{+}$ described before
$a=\frac{x_{1}^{+}}{2}$. We may assume that $x_{1}^{+}$ does not converge to
$0$. Let $x_{0}(\varepsilon_{m})$ be the first point where $v_{\varepsilon
_{m},n}^{+}(x_{0})=\varepsilon_{m}$ or $x_{0}=a$ if such a point does not
exist. We claim that $x_{0}(\varepsilon_{m})\rightarrow0$, as $\varepsilon
_{m}\rightarrow0$. It is clear that ${\partial^{2}v_{\varepsilon_{m},n}^{+}%
}/{\partial x^{2}}=-f_{\varepsilon_{m}}(v_{\varepsilon_{m},n}^{+})<0$ in
$(0,x_{1}^{+})$, and then
\begin{equation}
\displaystyle\frac{v_{\varepsilon_{m},n}^{+}(x_{0})}{x_{0}}x\leq
v_{\varepsilon_{m},n}^{+}(x)\leq\varepsilon_{m},\ \ \forall x\in\lbrack
0,x_{0}],\label{concavidadv}%
\end{equation}
by concavity. Hence, integrating first on $(s,a)$ and then on $(0,x)$ with
$x\leq x_{0},$ we have
\begin{equation}%
\begin{array}
[c]{l}%
\displaystyle\frac{d}{dx}v_{\varepsilon_{m},n}^{+}(s)=\int_{s}^{a}%
f_{\varepsilon_{m}}(v_{\varepsilon_{m},n}^{+}(\tau))d\tau,
\end{array}
\label{integrales}%
\end{equation}%
\[%
\begin{array}
[c]{l}%
v_{\varepsilon_{m},n}^{+}(x)=\displaystyle\int_{0}^{x}\int_{x_{0}}%
^{a}f_{\varepsilon_{m}}(v_{\varepsilon_{m},n}^{+}(\tau))d\tau
ds+\displaystyle\int_{0}^{x}\int_{s}^{x_{0}}f_{\varepsilon_{m}}(v_{\varepsilon
_{m},n}^{+}(\tau))d\tau ds.
\end{array}
\]

Since $f_{\varepsilon}(u)$ is concave, we have that $f_{\varepsilon}(u)/u\geq
f_{\varepsilon}(\varepsilon)/\varepsilon$, $\forall\ 0<u\leq\varepsilon$.
Moreover, by assumption ($\widetilde{A5}$) of $f_{\varepsilon}$ we get
$f_{\varepsilon}(u)\geq\frac{1-\varepsilon}{\varepsilon}u,\forall
0<u\leq\varepsilon$. Hence, using (\ref{concavidadv}) we have
\[%
\begin{array}
[c]{l}%
v_{\varepsilon_{m},n}^{+}(x)\geq\displaystyle\int_{0}^{x}\int_{s}^{x_{0}%
}\displaystyle\frac{1-\varepsilon_{m}}{\varepsilon_{m}}v_{\varepsilon_{m}%
,n}^{+}(\tau)d\tau ds\geq\displaystyle\frac{1-\varepsilon_{m}}{\varepsilon
_{m}}\displaystyle\frac{v_{\varepsilon_{m},n}^{+}(x_{0})}{x_{0}}%
\displaystyle\int_{0}^{x}\int_{s}^{x_{0}}\tau d\tau ds.
\end{array}
\]
Thus,
\[
1\geq\displaystyle\frac{1-\varepsilon_{m}}{\varepsilon_{m}}\left(
\frac{xx_{0}}{2}-\frac{x^{3}}{6x_{0}}\right)  ,
\]
so it follows that $x_{0}\rightarrow0,$ as $\varepsilon_{m}\rightarrow0.$

Let
$\delta_{1}<0<\delta_{2}$ be such that $x_{0}(\varepsilon_{m})\leq\delta
_{1}<\delta_{2}\leq a(\varepsilon_{m}).$ Since $v_{\varepsilon_{m},n}%
^{+}(x)\geq\varepsilon_{m}\ \forall x\in\lbrack x_{0},a],$ if we intregate
(\ref{integrales}) over $(\delta,x)$ with $\delta_{1}<x\leq\delta_{2}$, we
have
\[
v_{\varepsilon_{m},n}^{+}(x)-v_{\varepsilon_{m},n}^{+}(\delta
)=\displaystyle\int_{\delta_{1}}^{x}\int_{s}^{a}f(v_{\varepsilon_{m},n}%
^{+}(\tau))d\tau ds\geq(1-\varepsilon_{m})\displaystyle\int_{\delta_{1}}%
^{x}\int_{s}^{a}d\tau ds,
\]
which implies a contradiction if $v_{\varepsilon_{m},n}^{+}\rightarrow0$ in
$C([0,1]).$

The proof is similar for $v_{\varepsilon_{m},n}^{-}.$
\end{proof}

\begin{lemma}
\label{lemaconvergencia} $v_{\varepsilon_{m},k}^{+}$ (resp. $v_{\varepsilon
_{m},k}^{-})$ converges to $v_{k}^{+}$ (resp. $v_{k}^{-}$) in $H_{0}%
^{1}(\Omega)$ as $m\rightarrow\infty$ for any $k\geq1$.
\end{lemma}

\begin{proof}
It is easy to see that $v_{\varepsilon_{m}k}^{+}$ is bounded in $H^{2}%
(\Omega)\cap H_{0}^{1}(\Omega),$ so $v_{\varepsilon_{m}k}^{+}\rightarrow v$
strongly in $H_{0}^{1}(\Omega)$ and $C([0,1]).$ The proof will be finished if
we prove that $v=v_{k}^{+}$. It is clear that the functions $g_{\varepsilon
_{m}}=f_{\varepsilon_{n}}(v_{\varepsilon_{m}k}^{+})$ are bounded in
$L^{\infty}(0,1)$.

Passing to a subsequence we can then assume that $g_{\varepsilon_{n}}$
converges to some $g$ weakly in $L^{2}(0,1).$ It is clear that $-(\partial
^{2}v/\partial x^{2})=g$ and $v$ is a fixed point if we prove the inclusion
$g(x)\in H_{0}(v(x))$ for a.e. $x\in(0,1)$. By Masur's theorem \cite[p.120]%
{yosida} there exist $z_{m}\in V_{m}=conv(\cup_{k\geq m}^{\infty
}g_{\varepsilon_{k}})$ such that $z_{m}\rightarrow g,$ as $m\rightarrow
\infty,$ strongly in $L^{2}(0,1).$ Taking a subsequence we have $z_{m}%
(x)\rightarrow g(x),$ a.e. in $(0,1)$. Since $z_{m}\in V_{m},$ we get
$z_{m}=\sum_{i=1}^{N_{m}}\lambda_{i}g_{\varepsilon_{k_{i}}},$ where
$\lambda_{i}\in\lbrack0,1],\sum_{i=1}^{N_{m}}\lambda_{i}=1$ and $k_{i}\geq m,$
for all $i.$

Now (\ref{adecuacion}) implies that $|g_{\varepsilon_{k}}(x)-H_{0}%
(v(x))|\rightarrow0,$ as $k\rightarrow\infty$, for a.e. $x.$ Indeed, if
$v(x)=0,$ then $g_{\varepsilon_{k}}(x)\in\lbrack-1,1]=H_{0}(v(x)).$ If
$v(x)>0,$ then $|g_{\varepsilon_{k}}(x)-H_{0}(v(x))|=|f_{\varepsilon_{k}%
}(v_{\varepsilon_{k}}(x))-1|\rightarrow0,$ as $k\rightarrow\infty.$ If
$v(x)<0$, we apply a similar argument. 

Thus, for any $\delta>0$ and a.e. $x$
there exists $m(x,\delta)$ such that $g_{\varepsilon_{k}}(x)\subset\lbrack
a(x)-\delta,b(x)+\delta],$ for all $k\geq m,$ where $[a(x),b(x)]=H_{0}(v(x)).$
Hence, $z_{m}(x)\subset\lbrack a(x)-\delta,b(x)+\delta]$, as well. Passing to
the limit we obtain $g(x)\in\lbrack a(x),b(x)],$ a.e. on $(0,1).$

To conclude the proof, we have to prove that $v=v_{k}^{+}.$ By Lemma
\ref{lemaconvergenciaA0} $v\not =0.$ Hence, as $v_{\varepsilon_{m}k}^{+}(x)>0$
for all $x\in(0,x_{1}^{+}(\varepsilon_{m}))$, $v=v_{n}^{+}$ for some
$n\in\mathbb{N}$. Since $v_{n}^{+}$ has $n+1$ zeros, the convergence
$v_{\varepsilon_{m}k}^{+}\rightarrow v_{n}^{+}$ implies that $v_{\varepsilon
_{m}k}^{+}$ has $n+1$ zeros for $m\geq N$. But $v_{\varepsilon_{m}k}^{+}$
possesses $k+1$ zeros. Thus, $k=m$. 

For the sequence $v_{\varepsilon_{m}k}%
^{-}$ the proof is analogous.
\end{proof}

Once we have described the preliminary properties, we are now ready to check
that (\ref{aproximaciones}) satisfies the conditions given in Theorem
\ref{teoremaprincipal} for certain families $\mathcal{M}_{\varepsilon}$. We
recall that \cite[Theorem 10]{valerokapustyan} guarantees the existence of the
global compact invariant attractors $\mathcal{A}_{\varepsilon}$, where each
$\mathcal{A}_{\varepsilon}$ is the union of all bounded complete trajectories.

Let us check assumptions (H1)-(H5) of Theorem \ref{teoremaprincipal}.

As we have seen before, condition (H2) follows from Remark
\ref{condicion2delmaintheorem}. Therefore, we prove now condition (H1).

Multiplying the equation in (\ref{aproximaciones}) by $u$, we obtain
\begin{align}
\displaystyle \frac{1}{2}\frac{d}{dt}\|u\|^{2}_{L^{2}}+\|u\|^{2}_{H^{1}_{0}}
&  \leq\int_{\Omega}|u|dx\nonumber\\
&  \leq\frac{1}{2} \|u\|^{2}_{H^{1}_{0}}+C, \label{primerades}%
\end{align}
where we have used Poincar\'e's inequality, denoting $\lambda_{1}$ the first
eigenvalue of the operator $- \Delta$ in $H_{0}^{1}(\Omega)$.

Hence, we have
\[
\frac{d}{dt}\|u\|^{2}_{L^{2}}\leq- \lambda_{1}\|u\|^{2}_{L^{2}} + K.
\]
Gronwall's lemma gives
\begin{equation}
\label{cotal2}\|u(t)\|^{2}_{L^{2}} \leq e^{-\lambda_{1}t}\|u(0)\|^{2}_{L^{2}%
}+\frac{K}{\lambda_{1}}, \ \ t \geq0.
\end{equation}

Integrating (\ref{primerades}) over $(t,t+r)$ with $r>0$ we have
\[
\|u(t+r)\|_{L^{2}}^{2}+\int_{t}^{t+r}\|u\|^{2}_{H^{1}_{0}} ds\leq
\|u(t)\|_{L^{2}}^{2}+rK
\]
Then by (\ref{cotal2}),
\begin{equation}
\label{cotaenh01}\int_{t}^{t+r}\|u\|^{2}_{H^{1}_{0}} ds \leq\|u(0)\|_{L^{2}%
}^{2}e^{-\lambda_{1}t}+\left(  \frac{1}{\lambda_{1}}+r \right)  K.
\end{equation}

On the other hand, multiplying (\ref{aproximaciones}) by $- \Delta u$ and
using Young's inequality we obtain
\begin{equation}
\label{25b}\frac{d}{dt} \|u\|^{2}_{H^{1}_{0}}+2\|\Delta u\|_{L^{2}}^{2}%
\leq\|f_{\varepsilon}(u)\|_{L^{2}}^{2}+\|\Delta u\|_{L^{2}}^{2}%
\end{equation}

Since $f_{\varepsilon}(u(\cdot)) \in L^{2}(0,T; L^{2}(\Omega)), \forall T >0,
$ we obtain by \cite[p.189]{barbu} that
\[
u \in L^{\infty}(\eta,T; H_{0}^{1}(\Omega)),
\]
\[
\displaystyle \frac{du}{dt} \in L^{2}(\eta,T; L^{2}(\Omega)), \ \ \ \forall\ 0
< \eta<T.
\]

This regularity guarantees that the equality
\begin{equation}
\label{derivadadelcuadrado}\displaystyle \frac{1}{2} \frac{d}{dt}
\|u\|^{2}_{H^{1}_{0}}= \langle\frac{du}{dt}, -\Delta u \rangle, \text{ for
a.e. } t,
\end{equation}
is correct \cite[p.102]{sellyou}. Then
\[
\frac{d}{dt}\|u\|^{2}_{H^{1}_{0}}\leq\overline{K}+\|u\|^{2}_{H^{1}_{0}}.
\]

We apply the uniform Gronwall lemma \cite[p. 91]{temmam} with
$y(s)=\|u(s)\|^{2}_{H^{1}_{0}}$, $g(s)=1$ and $w(s)=\overline{K}$. Also, using
(\ref{cotaenh01}) we obtain
\begin{equation}
\label{cotatr}\|u(t+r)\|^{2}_{H^{1}_{0}} \leq\left(  \frac{\|u(0)\|_{L^{2}%
}^{2}e^{-\lambda_{1}t}+ (\frac{1}{\lambda_{1}}+r)K}{r}+\overline{K} \right)
e^{r}%
\end{equation}

It follows from (\ref{cotal2}) that $\| y\|_{L^{2}} \leq\frac{K}{\lambda_{1}}$
for any $y \in\mathcal{A_{\varepsilon}}$, $0 <\varepsilon\leq1$. Hence,
$\cup_{0< \varepsilon\leq1} \mathcal{A}_{\varepsilon}$ is bounded in
$L^{2}(\Omega)$. Since $\mathcal{A_{\varepsilon}} \subset G_{\varepsilon}(t,
\mathcal{A_{\varepsilon}})$ for any $t \geq0$, for any $y\in\mathcal{A}%
_{\varepsilon}$ there exists $z \in\mathcal{A}_{\varepsilon} $ such that $y
\in G_{\varepsilon}(1,z)$. Then using (\ref{cotatr}) with $r=1 $ and $t=0$ we
obtain that
\[
\|y\|^{2}_{H^{1}_{0}} \leq\left(  \|z\|_{L^{2}}^{2}+ \left(
\displaystyle \frac{1}{\lambda_{1}}+1 \right)  K+\overline{K} \right)  e,
\]
so $\cup_{0< \varepsilon\leq1} \mathcal{A}_{\varepsilon}$ is bounded in
$H^{1}_{0}(\Omega)$. The compact embedding $H_{0}^{1} (\Omega) \subset
L^{2}(\Omega)$ implies that $\cup_{0< \varepsilon\leq1} \mathcal{A}%
_{\varepsilon}$ is relatively compact in $L^{2}(\Omega)$. As the global
attractor $A_{0}$ of the differential inclusion (\ref{inclusion1}) is compact,
the set $\overline{\cup_{0\leq\varepsilon\leq1} \mathcal{A}_{\varepsilon}}$ is
compact in $L^{2}(\Omega)$.

In order to establish that (\ref{aproximaciones}) satisfies the rest of
conditions given in Theorem \ref{teoremaprincipal}, we need to proof two
previous results related to the convergence of solutions of the approximations
and the connections between fixed points.

\begin{theorem}
\label{convergenciatrayectorias} If $u_{\varepsilon_{n}0}\rightarrow u_{0}$ in
$L^{2}(\Omega)$ as $\varepsilon_{n}\rightarrow0,$ then for any sequence of
solutions of (\ref{aproximaciones}) $u_{\varepsilon_{n}}(\cdot)$ with
$u_{\varepsilon_{n}}(0)=u_{\varepsilon_{n}0}$ there exists a subsequence of
$\varepsilon_{n}$ such that $u_{\varepsilon_{n}}$ converges to some strong
solution $u$ of (\ref{inclusion1}) in the space $C([0,T],L^{2}(\Omega))$, for
any $T>0.$
\end{theorem}

\begin{proof}
We define $g_{n}(t)=f_{\varepsilon_{n}}(u_{\varepsilon_{n}}(t))$ and $u_{n}(t)=u_{\varepsilon_{n}}(t)$. From (\ref{cotal2})  we have that $\|u_{n}(t)\|_{L^{2}} \leq C_{0},$ for all $t \geq 0,$ so that $\|g_{n}(t)\|_{L^{2}}\leq C_{1},$ for a.e. $t \geq 0.$ Hence, there exists a subsequence such that $u_{n}\rightarrow u$ weakly in $L^{2}(0,T; L^{2}(\Omega)).$
It follows from (\ref{25b}) and $\|g_{n}(t) \|_{L^{2}} \leq C_{1}$ that $\int_{r}^{T} \|\Delta u\|_{L^{2}}^{2}ds \leq C_{1}^{2}(T-r)+\|u_{n}(r)\|^{2}_{H_{0}^{1}}.$ Using (\ref{cotatr}) we obtain that $\int_{r}^{T} \|\Delta u_{n}\|_{L^{2}}^{2}ds \leq C(r).$ Hence, $\frac{du_{n}}{dt}$ is bounded in $L^{2}(r,T;L^{2}(\Omega))$ for any $0<r<T,$ so passing to a subsequence $\frac{du_{n}}{dt} \rightarrow \frac{du}{dt}$ weakly in $L^{2}(r,T;L^{2}(\Omega)).$

Moreover, Ascoli-Arzel\`a theorem implies that for any fixed $r>0$ we have $u_{n} \rightarrow u$ in $C([r,T],L^{2}(\Omega))$ and $u$ is absolutely continuous on $[r,T].$

Also, $g_{n}$ converges to some $g \in L^{\infty}(0,T;L^{2}(\Omega))$ weakly star in $L^{\infty}(0,T;L^2(\Omega))$ and weakly in $L^{2}(0,T;L^{2}(\Omega)).$ On the other hand, since $-\Delta u_{n}=-\frac{du_{n}}{dt}+g_{n}$, $-\Delta u_{n}$ converges to $l(t)=-(\frac{du}{dt})+g$ weakly in $L^{2}(r,T;L^{2}(\Omega)).$ Hence, we find at once that $u$ satisfies $$\frac{du}{dt}-\Delta u(t)=g(t), \textrm{ a.e. on } (0,T).$$

The result will be established if we proof that $u(\cdot)$ is a strong solution of (\ref{inclusion1}). Now, we show that $g(t) \in H_{0}(u(t)),$ a.e. in $(0,T).$ For this, we shall prove first that for a.e. $x\in \Omega$ and $s\in (0,T)$
$$|g_{n}(s,x)-H_{0}(u(s,x)) | \rightarrow 0, \textrm{ as } n \rightarrow \infty.$$

Indeed, if $u(s,x)=0,$ then $g_{n}(s,x)=f_{\varepsilon_{n}}(u_{n}(s,x))=0 \in [-1,1]=H_{0}(u(s,x)),$ for all $n$, so that the result is evident. If $u(s,x)<0,$ then $$|g_{n}(s,x)-H_{0}(u(s,x)) |= |f_{\varepsilon_{n}}(u_{n}(s,x))+1 | \rightarrow 0, \textrm{ as } n \rightarrow \infty.$$
Finally, if $u(s,x)>0,$ then $$|g_{n}(s,x)-f_{0}(u(s,x)) |= |f_{\varepsilon_{n}}(u_{n}(s,x))-1 | \rightarrow 0, \textrm{ as } n \rightarrow \infty.$$

Now, by \cite[Proposition 1.1]{tolstonogov} we  have that for a.e. $t \in (0,T)$
$$g(t) \in \bigcap_{n\geq 0}\overline{co} \bigcup_{k\geq n}g_{k}(t).$$
Then $g(t)=\displaystyle \lim_{n\rightarrow \infty} y_{n}(t)$ strongly in $L^{2}(\Omega),$ where $$y_{n}(t)=\displaystyle \sum_{i=1}^{M} \lambda_{i}g_{k_{i}}(t), \displaystyle \sum_{i=1}^{M}\lambda_{i}=1, k_{i} \geq n.$$

We note that for any $t \in [0,T]$ and a.e. $x \in \Omega$ we can find $n(\varepsilon,x,t)$ such that if $k\geq n,$ then $|g_{k}(t,x)-H_{0}(u(t,x)) | \leq \varepsilon.$ Therefore,
$$|y_{n}(t,x)-H_{0}(u(t,x)) | \leq  \displaystyle \sum_{i=1}^{M} \lambda_{i} |g_{k_{i}}(t,x)-H_{0}(u(t,x))| \leq \varepsilon.  $$
Hence, since we can assume that for a.e. $(t,x) \in (0,T) \times \Omega, y_{n}(t,x) \rightarrow g(t,x),$ it follows that $g(t,x) \in H_{0}(u(t,x)).$

It remains to check that $u$ is continuous as $t \rightarrow 0^{+}.$ Let $\hat{u} $ be the unique solution of \begin{equation*}
\left\{
\begin{array}
[c]{c}
\displaystyle \frac{du}{dt}-\Delta u=0,\\
u|_{\partial \Omega}=0,\\
u(0)=u_{0},
\end{array}
\right.
\end{equation*}
and let $v_{n}(t)=u_{n}(t)-\hat{u}(t).$ Multiplying by $v_{n}$ the equation $$\frac{dv_{n}}{dt}-\Delta v_{n}=f_{\varepsilon_{n}}(u_{n}),$$ we obtain
\begin{eqnarray*}
\frac{1}{2} \frac{d}{dt} \|v_{n} \|_{L^{2}}^{2}+\|v_{n} \|^{2}_{H^{1}_{0}} \leq (f_{\varepsilon_{n}}(u_{n}(t)),v_{n})  \\
\leq \frac{1}{2}\|f_{\varepsilon_{n}}(u_{n}) \|_{L^{2}}^{2}+\frac{1}{2} \|v_{n} \|_{L^{2}}^{2} \leq
C \|v_{n} \|_{L^{2}}^{2} \leq K,
\end{eqnarray*}
so that $$\|v_{n}(t) \|_{L^{2}}^{2} \leq \|v_{n}(0) \|_{L^{2}}^{2}+Kt. $$

Hence, $\|u(t)-\hat{u}(t) \|_{L^{2}}^{2}= \lim_{n\rightarrow \infty} \|v_{n}(t) \|_{L^{2}}^{2} \leq Kt,$ for $t>0$, and $$\|u(t)-u_{0} \|_{L^{2}} \leq \|u(t)-\hat{u}(t) \|_{L^{2}}+ \|\hat{u}(t)-u_{0} \| _{L^{2}}< \delta, $$ as soon as $t<\varepsilon (\delta).$
Therefore, $u(\cdot)$ is a strong solution.

Finally, if $t_{n}\rightarrow0$, then
\begin{align*}
\left\Vert u_{n}(t_{n})-u_{0}\right\Vert _{L^{2}}  & \leq\left\Vert
v_{n}(t_{n})\right\Vert _{L^{2}}+\left\Vert \widehat{u}(t_{n})-u_{0}%
\right\Vert _{L^{2}}\\
& \leq\sqrt{\left\Vert v_{n}(0)\right\Vert _{L^{2}}^{2}+Kt_{n}}+\left\Vert
\widehat{u}(t_{n})-u_{0}\right\Vert _{L^{2}}\rightarrow0.
\end{align*}
Hence, $u_{n}\rightarrow u$ in $C([0,T],L^{2}(\Omega))$.
\end{proof}

As a consequence of the last theorem, condition (H4) follows.

\begin{remark}
\label{fantasma} Let be $u_{\varepsilon_{n}}(\cdot)$ a bounded complete
trajectory of (\ref{aproximaciones}). Fix $T>0.$ Since $\bigcup_{0<\varepsilon
\leq\varepsilon_{0}}\mathcal{A_{\varepsilon}}$ is precompact in $L^{2}%
(\Omega)$, $u_{\varepsilon_{n}}(-T)\rightarrow y$ in $L^{2}$ up to a
subsequence. Theorem \ref{convergenciatrayectorias} implies that
$u_{\varepsilon_{n}}$ converges in $C([0,T],L^{2}(\Omega))$ to some solution
$u$ of (\ref{inclusion1}). If we choose successive subsequences for
$-2T,-3T,\ldots,$ and apply the standard diagonal procedure, we obtain that a
subsequence $u_{\varepsilon_{n}}$ converges to a complete trajectory $u$ of
(\ref{inclusion1}) in $C([-T,T],L^{2}(\Omega))$ for any $T>0.$
\end{remark}

Now, we need to prove a previous lemma to obtain the convergence of solutions
of the approximations in the space $C([0,T],H_{0}^{1}).$

\begin{lemma}
\label{precompacidadH10}Any sequence $\xi_{n}\in A_{\varepsilon_{n}}$ with
$\varepsilon_{n}\rightarrow0$ is relatively compact in $H_{0}^{1}(\Omega).$
\end{lemma}

\begin{proof}
For fixed $t_{0}\in\mathbb{R}$ there exists a bounded complete trajectory
$\psi_{\varepsilon_{n}}$ of (\ref{aproximaciones}) with $\psi_{\varepsilon
_{n}}(t_{0})=\xi_{n}.$
Denote $u_{n}($\textperiodcentered$)=\psi_{\varepsilon_{n}}(t_{0}%
+$\textperiodcentered$)$ and choose some $T>0$. Then $\xi_{n}=u_{n}(T)$,
$u_{n}(0)=\psi_{\varepsilon_{n}}(t_{0}-T)$. In view of Remark \ref{fantasma}
up to a subsequence $u_{n}\rightarrow u$ in $C([0,T],L^{2}(\Omega))$, where
$u$ is a strong solution of (\ref{inclusion1}). On top of that, by
(\ref{cotatr}) and the argument in the proof of Theorem
\ref{convergenciatrayectorias} we obtain that for $r>0$,
\begin{align*}
u_{n}  & \rightarrow u\text{ weakly star in }L^{\infty}(r,T;H_{0}^{1}%
(\Omega)),\\
\frac{du_{n}}{dt}  & \rightarrow\frac{du}{dt}\text{ weakly in }L^{2}%
(r,T;L^{2}(\Omega)),\\
u_{n}  & \rightarrow u\text{ weakly in }L^{2}(r,T;H^{2}(\Omega)).
\end{align*}
Therefore, by the Compactness Theorem \cite[p.58]{lions} we have
\begin{align*}
u_{n}  & \rightarrow u\text{ strongly in }L^{2}(r,T,H_{0}^{1}(\Omega)),\\
u_{n}(t)  & \rightarrow u\left(  t\right)  \text{ in }H_{0}^{1}(\Omega)\text{
for a.a. }t\in(r,T).
\end{align*}
In addition, by standard results \cite[p.102]{sellyou} we have that
$u_{n},u\in C([r,T],H_{0}^{1}(\Omega)).$

Multiplying (\ref{aproximaciones}) by $\frac{du_{n}}{dt}$ and using
(\ref{derivadadelcuadrado}), we obtain
\[
\displaystyle\left\Vert \frac{du_{n}}{dt}\right\Vert _{L^{2}}^{2}+\frac{d}%
{dt}\Vert u_{n}\Vert_{H_{0}^{1}}^{2}\leq\Vert f_{\varepsilon}(u_{n}%
)\Vert_{L^{2}}^{2}.
\]
Thus,
\[
\Vert u_{n}(t)\Vert_{H_{0}^{1}}^{2}\leq\Vert u_{n}(s)\Vert_{H_{0}^{1}}%
^{2}+C(t-s),\ C>0,\ t\geq s\geq r.
\]
The same inequality is valid for the limit function $u(\cdot).$ Hence, the
functions $J_{n}(t)=\Vert u_{n}(t)\Vert_{H_{0}^{1}}^{2}-Ct,J(t)=\Vert
u(t)\Vert_{H_{0}^{1}}^{2}-Ct$, are continuous and non-increasing in $[r,T].$
Moreover, $J_{n}(t)\rightarrow J(t)$ for a.e. $t\in(r,T).$ Take $r<t_{m}<T$
such that $t_{m}\rightarrow T$ and $J_{n}(t_{m})\rightarrow J(t_{m})$ for all
$m.$ Then
\[
J_{n}(T)-J(T)\leq J_{n}(t_{m})-J(T)\leq|J_{n}(t_{m})-J(t_{m})|+|J(t_{m}%
)-J(T)|.
\]
For any $\varepsilon>0$ there exist $m(\varepsilon)$ and $N(\varepsilon)$ such
that $J_{n}(T)-J(T)\leq\varepsilon$ if $n\geq N.$ Then $\limsup J_{n}(T)\leq
J(T),$ so $\limsup\Vert u_{n}(T)\Vert_{H_{0}^{1}}^{2}\leq\Vert u(T)\Vert
_{H_{0}^{1}}^{2}.$ As $u_{n}(T)\rightarrow u(T)$ weakly in $H_{0}^{1}$ implies
$\liminf\Vert u_{n}(T)\Vert_{H_{0}^{1}}^{2}\geq\Vert u(T)\Vert_{H_{0}^{1}}%
^{2},$ we obtain
\[
\Vert u_{n}(T)\Vert_{H_{0}^{1}}^{2}\rightarrow\Vert u(T)\Vert_{H_{0}^{1}}^{2},
\]
so that $u_{n}(T)\rightarrow u(T)$ strongly in $H_{0}^{1}(\Omega).$ Hence, the
result follows.
\end{proof}

\begin{corollary}
\label{corolario} If $u_{\varepsilon0}\rightarrow u_{0}$ in $L^{2}(\Omega)$,
where $u_{\varepsilon0}\in\mathcal{A}_{\varepsilon}$, $u_{0}\in\mathcal{A}%
_{0},$ then for any $T>0$ there exists a subsequence $\varepsilon_{n}$ such
that $u_{\varepsilon_{n}}$ converges to some strong solution $u$ of
(\ref{inclusion1}) in $C([0,T],H_{0}^{1}(\Omega)).$
\end{corollary}

We know from Theorem \ref{convergenciatrayectorias} that there exists a
subsequence such that $u_{\varepsilon_{n}}$ converges to some strong solution
$u$ of (\ref{inclusion1}) in $C([0,T],L^{2}(\Omega))$. Then the statement
follows from the invariance of $\mathcal{A}_{\varepsilon}$ and Lemma
\ref{precompacidadH10}.

\begin{remark}
\label{convergenciaaproxfinal} Let be $u_{\varepsilon_{n}}(\cdot)$ a complete
trajectory of (\ref{aproximaciones}). Fix $T>0$. By Lemma
\ref{precompacidadH10} $u_{\varepsilon_{n}}(-T)\rightarrow y$ in $H_{0}%
^{1}(\Omega)$ up to a subsequence. Corollary \ref{corolario} implies then that
$u_{\varepsilon_{n}}$ converges in $C([0,T],H_{0}^{1}(\Omega))$ to some
solution $u$ of (\ref{inclusion1}). If we choose successive subsequences for
$-2T,-3T\ldots$ and apply the standard diagonal procedure we obtain that a
subsequence $u_{\varepsilon_{n}}$ converges to a complete trajectory $u$ of
(\ref{inclusion1}) in $C([-T,T],H_{0}^{1}(\Omega))$ for any $T>0.$
\end{remark}

\begin{lemma}
\label{lemaatractores} $dist_{H_{0}^{1}}(\mathcal{A}_{\varepsilon}%
,\mathcal{A}_{0})\rightarrow0$, as $\varepsilon\rightarrow0.$
\end{lemma}

\begin{proof}
By contradiction let there exist $\delta>0$ and a sequence $y_{\varepsilon
_{n}}\in\mathcal{A}_{\varepsilon_{n}}$ such that
\[
dist_{H_{0}^{1}}(y_{\varepsilon_{n}},\mathcal{A}_{0})>\delta.
\]
Hence, as $y_{\varepsilon_{n}}=u_{\varepsilon_{n}}(0),$ where $u_{\varepsilon
_{n}}$ is a bounded complete trajectory of problem (\ref{aproximaciones}),
using Remark \ref{convergenciaaproxfinal} we obtain that up to a sequence
$u_{\varepsilon_{n}}$ converges to a complete trajectory $u$ of the problem
(\ref{inclusion1}) in the spaces $C([-T,T],H_{0}^{1}(\Omega))$ for every
$T>0$. Since $\cup_{0<\varepsilon\leq1}\mathcal{A}_{\varepsilon}$ is bounded
in $L^{2}(\Omega)$ (in fact in $H_{0}^{1}(\Omega)),$ it is clear that $u$ is a
bounded complete trajectory of problem (\ref{inclusion1}). Thus,
$u(t)\in\mathcal{A}_{0}$ for any $t\in\mathbb{R}$. We infer then that
\[
y_{\varepsilon_{n}}=u_{\varepsilon_{n}}(0)\rightarrow u(0)\in\mathcal{A}_{0},
\]
which is a contradiction.
\end{proof}

\begin{remark}
\label{remarkatractor} Let be $u_{\varepsilon_{n}}(\cdot)$ a complete
trajectory of (\ref{aproximaciones}). By Remark \ref{convergenciaaproxfinal}
we have that $u_{\varepsilon_{n}}(\cdot)$ converges to a complete trajectory
$u$ of (\ref{inclusion1}) in $C([-T,T],H^{1}_{0}(\Omega))$ for any $T>0.$
Since $u_{\varepsilon_{n}}(t) \in\mathcal{A}_{\varepsilon_{n}}$, by Lemma
\ref{lemaatractores} $u(t)\in\mathcal{A}_{0}$ for all $t$. Therefore, $u(t)$
is a bounded complete trajectory contained in the global attractor.
\end{remark}

We choose some $\delta>0$ such that
\[
\mathcal{O}_{\delta}(\Xi_{i}^{0}) \cap\mathcal{O}_{\delta}(\Xi_{j}^{0})=
\emptyset\text{ if } i \ne j
\]
and $\Xi_{i}^{0}$ are maximal weakly invariant.

For problem (\ref{aproximaciones}) let us define the sets
\[
M_{i}^{\varepsilon}=\{v_{\varepsilon,i}^{+}, v_{\varepsilon,i}^{-}\} \text{
for } 1\leq i <n,
\]
\[
Z_{n}^{\varepsilon}= \left(  \displaystyle{\cup_{k \geq n}} \{v_{\varepsilon
,k}^{\pm}\} \right)  \displaystyle{\cup} \{0 \},
\]

\[
M_{n}^{\varepsilon}= \left\{
\begin{array}
[c]{c}%
y: \exists\psi\in\mathbb{K} ^{\varepsilon} \text{ such that
(\ref{trayectoriasconectando}) holds} \text{ with } z_{j} \in Z_{n}%
^{\varepsilon},\\
j=1,2 \text{ and } y=\psi(t) \text{ for some } t \in\mathbb{R}%
\end{array}
\right\}
\]
where $\mathbb{K}^{\varepsilon}$ is the set of all bounded complete
trajectories of (\ref{aproximaciones}).

In view of Lemma \ref{lemaconvergencia} we have
\[
dist_{H^{1}_{0}}(M_{i}^{\varepsilon}, \Xi_{i}^{0}) \rightarrow0, \text{ as }
\varepsilon\rightarrow0, \ \ 1 \leq i < n
\]

\begin{lemma}
\label{ConvergMorseSetN}$dist_{H_{0}^{1}}(M_{n}^{\varepsilon},\Xi_{n}%
^{0})\rightarrow0,\text{ as }\varepsilon\rightarrow0.$
\end{lemma}

\begin{proof}
Suppose the opposite, that is, there exists $\delta>0$ and a sequence
$y_{\varepsilon_{k}}$ such that
\begin{equation}
dist_{H_{0}^{1}}(y_{\varepsilon_{k}},\Xi_{n}^{0})>\delta\text{ for all
}k.\label{contradiccion}%
\end{equation}
Let $\xi_{\varepsilon_{k}}$ be a sequence of bounded complete trajectories of
problem (\ref{aproximaciones}) such that
\[
\xi_{\varepsilon_{k}}(t)\rightarrow z_{-1}^{k}\text{ as }t\rightarrow-\infty,
\]%
\[
\xi_{\varepsilon_{k}}(t)\rightarrow z_{0}^{k}\text{ as }t\rightarrow\infty,
\]
where $z_{-1}^{k},z_{0}^{k}\in Z_{n}^{\varepsilon_{k}}$. By Lemma
\ref{lemaconvergencia}, passing to a subsequence we have that
\[
z_{i}^{k}\rightarrow z_{i}\in Z_{n},i=-1,0.
\]

By Remark \ref{convergenciaaproxfinal} we obtain that up to a subsequence
$\xi_{\varepsilon_{k}}$ converges to a complete trajectory $\psi_{0}$ of
problem (\ref{inclusion1}) in the spaces $C([-T,T],H_{0}^{1}(\Omega))$ for
every $T>0$, so $y_{\varepsilon_{k}}\rightarrow\psi_{0}(0)$ in $H_{0}%
^{1}(\Omega).$ We know that there exists two fixed points of problem
(\ref{inclusion1}), denoted by $\overline{z}_{-1},\overline{z}_{0}$ such that
\[
E(\overline{z}_{-1})>E(\overline{z}_{0}),
\]%
\[
\psi_{0}(t)\rightarrow\overline{z}_{-1}\text{ as }t\rightarrow-\infty,
\]%
\[
\psi_{0}(t)\rightarrow\overline{z}_{0}\text{ as }t\rightarrow\infty.
\]

If $\overline{z}_{0}=z_{0},$ then $\overline{z}_{-1},\overline{z}_{0}\in
Z_{n},$ which means that $\psi_{0}(0)\in\Xi_{n}^{0}.$ This would imply a
contradiction with (\ref{contradiccion}). Therefore, we assume that
$\overline{z}_{0}\not =z_{0}.$ Also, it is clear that $\overline{z}_{0}%
=v_{m}^{\pm}\not =0,$ for some $m\in\mathbb{N}.$ 

Let $r_{0}>0$ be such that $\mathcal{O}_{r_{0}}(\overline{z}_{0}%
)\cap\mathcal{O}_{r_{0}}(z_{0})\not =\emptyset$ and $\mathcal{O}_{2r_{0}%
}(z_{0})$ does not contain any other fixed point of problem (\ref{inclusion1}%
). The previous convergences imply that for each $r\leq r_{0}$ there exist a
moment of time $t_{r}$ and $k_{r}$ such that $\xi_{\varepsilon_{k}}(t_{r}%
)\in\mathcal{O}_{r}(\overline{z}_{0})$ for all $k\geq k_{r}.$ On the other
hand, since $\xi_{\varepsilon_{k}}(t)\rightarrow z_{0}^{k},$ as $t\rightarrow
\infty,$ and $z_{0}^{k}\rightarrow z_{0},$ there exists $t_{r}^{\prime}>t_{r}$
such that
\[
\xi_{\varepsilon_{k_{r}}}(t)\in\mathcal{O}_{r_{0}}(\overline{z}_{0})\text{ for
all }t\in\lbrack t_{r},t_{r}^{\prime}),
\]%
\[
\Vert\xi_{\varepsilon_{k_{r}}}(t_{r}^{\prime})-\overline{z}_{0}\Vert_{L^{2}%
}=r_{0}.
\]

Let us consider two cases: 1) $t_{r}^{\prime}-t_{r}\rightarrow\infty
;2)|t_{r}^{\prime}-t_{r}|\leq C.$ We begin with the first case. We define the
sequence of bounded complete trajectories of problem (\ref{aproximaciones})
given by
\[
\xi_{k_{r}}^{1}(t)=\xi_{\varepsilon_{k_{r}}}(t+t_{r}^{\prime}).
\]

By Remark \ref{convergenciaaproxfinal} we can extract a subsequence of this
sequence converging to a bounded complete trajectory $\psi_{1}$ of problem
(\ref{inclusion1}). Since $t_{r}^{\prime}-t_{r}\rightarrow\infty$, we obtain
that $\psi_{1}(t)\in\mathcal{O}_{r_{0}}(\overline{z}_{0})$ for all $t\leq0.$
Since $\mathcal{O}_{2r_{0}}(z_{0})$ does not contain any other fixed point of
problem (\ref{inclusion1}), it follows that $\psi_{1}(t)\rightarrow
\overline{z}_{0}$ as $t\rightarrow-\infty.$ But $\Vert\psi_{1}(0)-\overline
{z}_{0}\Vert_{L^{2}}=r_{0},$ so $\psi_{1}$ is not a fixed point. Therefore,
$\psi_{1}(t)\rightarrow\overline{z}_{1}$ as $t\rightarrow\infty,$ where
$\overline{z}_{1}$ is a fixed point such that $E(\overline{z}_{1}%
)<E(\overline{z}_{0}).$ 

In the second case we define the sequence
\[
\xi_{k_{r}}^{1}(t)=\xi_{\varepsilon_{k_{r}}}(t+t_{r}).
\]
Passing to a subsequence we have that
\[
\xi_{k_{r}}^{1}(0)\rightarrow\overline{z}_{0},
\]%
\[
t_{r}^{\prime}-t_{r}\rightarrow t^{\prime}.
\]
As $\xi_{k_{r}}^{1}$ converges to a solution $\xi^{1}$ of problem
(\ref{inclusion1}) uniformly in bounded subsets from $[0,\infty)$ such that
$\xi^{1}(0)=\overline{z}_{0}$, $\xi_{k_{r}}^{1}(t_{r}^{\prime}-t_{r}%
)\rightarrow\xi^{1}(t^{\prime})$, so that $\Vert\xi^{1}(t^{\prime}%
)-\overline{z}_{0}\Vert_{L^{2}}=r_{0}.$ We put
\[
\psi_{1}(t)=\left\{
\begin{array}
[c]{c}%
\overline{z}_{0}\text{ if }t\leq0,\\
\xi^{1}(t)\text{ if }t\geq0.
\end{array}
\right.
\]
Then $\psi_{1}$ is a bounded complete trajectory of problem (\ref{inclusion1})
such that $\psi_{1}(t)\rightarrow\overline{z}_{1}$ as $t\rightarrow\infty$,
where $\overline{z}_{1}$ is a fixed point satisfying $E(\overline{z}%
_{1})<E(\overline{z}_{0}).$ 

Now, if $\overline{z}_{1}=z_{0},$ then we have the
chain of connections
\[
\psi_{0}(t)\rightarrow\overline{z}_{-1}\text{ as }t\rightarrow-\infty,\psi
_{0}(t)\rightarrow\overline{z}_{0}\text{ as }t\rightarrow+\infty,
\]%
\[
\psi_{1}(t)\rightarrow\overline{z}_{0}\text{ as }t\rightarrow-\infty,\psi
_{1}(t)\rightarrow\overline{z}_{1}\text{ as }t\rightarrow+\infty,
\]
which implies that $\overline{z}_{-1},\overline{z}_{0},\overline{z}_{1}\in
Z_{n},$ an then $\psi_{0}(0)\in\Xi_{n}^{0}.$ This would imply a contradiction
with (\ref{contradiccion}). 

However, if $\overline{z}_{1}\not =\overline
{z}_{0},$ then we proceed in the same way and obtain a new connection from the
point $\overline{z}_{1}$ to another fixed point with less energy. Since the
number of fixed points with energy less than or equal to $E(\overline{z}_{0})$
is finite, we will finally obtain a chain of connections of the form
\begin{align*}
\psi_{0}(t) &  \rightarrow\overline{z}_{-1}\text{ as }t\rightarrow
-\infty,\ \psi_{0}(t)\rightarrow\overline{z}_{0}\text{ as }t\rightarrow
+\infty,\\
\psi_{1}(t) &  \rightarrow\overline{z}_{0}\text{ as }t\rightarrow
-\infty,\ \psi_{1}(t)\rightarrow\overline{z}_{1}\text{ as }t\rightarrow
+\infty,\\
&  \vdots\\
\psi_{n}(t) &  \rightarrow\overline{z}_{m-1}\text{ as }t\rightarrow
-\infty,\ \psi_{n}(t)\rightarrow\overline{z}_{m}=z_{0}\text{ as }%
t\rightarrow+\infty.
\end{align*}
And again, this implies a contradiction with (\ref{contradiccion}).
\end{proof}

These convergences imply the existence of $\varepsilon_{0}$ such that if
$\varepsilon\leq\varepsilon_{0}$, then
\[
M_{i}^{\varepsilon} \subset\mathcal{O}_{\delta}(\Xi_{i}^{0}) \text{ for any }
1 \leq i \leq n.
\]

Further, let
\[
\Xi_{i}^{\varepsilon}=\left\{  \
\begin{array}
[c]{c}%
y:\exists\psi\in\mathbb{K}^{\varepsilon}\text{ such that }\psi(0)=y\\
\text{ and }\psi(t)\in\mathcal{O}_{\delta}(\Xi_{i}^{0})\ \text{{for all }}%
t\in\mathbb{R}%
\end{array}
\right\}  .
\]

These sets are clearly maximal weakly invariant for $G_{\varepsilon}$ in
$\mathcal{O}_{\delta}(\Xi_{i}^{0})$, so condition (H5) is satisfied for
$V_{i}=\mathcal{O}_{\delta}(\Xi_{i}^{0})$. As a consequence of Lemmas
\ref{lemaconvergencia}, \ref{ConvergMorseSetN}, Remark \ref{fantasma} and the
definition of $\delta$ we have
\[
dist_{L^{2}}(\Xi_{i}^{\varepsilon},\Xi_{i}^{0})\rightarrow0,\text{ as
}\varepsilon\rightarrow0,\text{ for }1\leq i\leq n.
\]
Therefore, condition (H3) is satisfied.

We also get by Remark \ref{convergenciaaproxfinal} and the definition of
$\delta$ that
\[
dist_{H^{1}_{0}}(\Xi_{i}^{\varepsilon}, \Xi_{i}^{0}) \rightarrow0, \text{ as }
\varepsilon\rightarrow0, \text{ for } 1\leq i \leq n.
\]

Moreover, $\mathcal{M^{\varepsilon}}=\{\Xi_{1}^{\varepsilon}, \ldots, \Xi
_{n}^{\varepsilon}\}$ is a disjoint family of isolated weakly invariant sets.

Applying Theorem \ref{teoremaprincipal} we obtain the following result.

\begin{theorem}
There exists $\varepsilon_{1}>0$ such that for all $0<\varepsilon
\leq\varepsilon_{1}$ the multivalued semiflow $G_{\varepsilon}$ is dynamically
gradient with respect to the family $\mathcal{M}^{\varepsilon}$.
\end{theorem}

\subsection*{Acknowledgments}

This paper is dedicated to the memory of Professor Valery
Melnik, on the tenth anniversary of his passing away, with our deepest respect
and sorrow.


\begin{thebibliography}{99}                                                                                               %


\bibitem {aragaocosta2011}E. R. Arag\~ao-Costa, T. Caraballo, A. N. Carvalho
and J. A. Langa, Stability of gradient semigroups under perturbations,
\emph{Nonlinearity}, \textbf{24} (2011), 2099-2117.

\bibitem {aragaocosta2013}E. R. Arag\~ao-Costa, T. Caraballo, A. N. Carvalho
and J. A. Langa, Non-autonomous Morse-decomposition and Lyapunov functions for
gradient-like processes. \emph{Transactions of the American Mathematical
Society (10)}, \textbf{365} (2013), 5277-5312.

\bibitem {arrietavalero}J. M. Arrieta, A. Rodr\'iguez-Bernal and J. Valero,
Dynamics of a reaction-diffusion equation with a discontinuous nonlinearity,
\emph{International Journal of Bifurcation and Chaos (10)}, \textbf{16}
(2006), 2965-2984.

\bibitem {ball}J. M. Ball, Continuity properties and global attractors of
generalized semiflows and the Navier-Stokes equations, \emph{Journal of
Nonlinear Science}, \textbf{7} (1997), 475-502.

\bibitem {barbu}V. Barbu, {\itshape Nonlinear Semigroups and Differential
Equations in Banach Spaces}, Editura Academiei, Bucuresti, 1976.

\bibitem {caraballorubio}T. Caraballo, P. Mar\'in-Rubio and J. Robinson, A
comparison between two theories for multi-valued semiflows and their
asymptotic behaviour, \emph{Set-Valued Analysis}, \textbf{11} (2003), 297-322.

\bibitem {carvalholanga}A. N. Carvalho and J. A. Langa, An extension of the
concept of gradient semigroups which is stable under perturbation, \emph{J.
Differential Equations}, \textbf{246} (2009), 2646-2668.

\bibitem {chepvishik}V. V. Chepyzhov and M. I. Vishik, {\itshape Attractors
for Equations of Mathematical Physics}, Americal Mathematical Society,
Providence, 2002.

\bibitem {costavalero}H. B. da Costa and J. Valero, Morse decompositions and
Lyapunov functions for dynamically gradient multivalued semiflows,
\emph{Nonlinear Dyn.}, \textbf{84} (2016), 19-34.

\bibitem {gruber}P. Gruber, {\itshape Convex and Discrete Geometry}, Springer, 2007.

\bibitem {henry1981}D. Henry, {\itshape Geometric Theory of Semilinear
Parabolic Equations}, Springer, Berlin, 2007.

\bibitem {henry1985}D. Henry, Some infinite-dimensional Morse-Smale systems
defined by parabolic partial differential equations, \emph{J. Diff. Eqs.},
\textbf{59} (2007), 165-205.

\bibitem {kapustyankasyanov}O. V. Kapustyan, P. O. Kasyanov and J. Valero,
Structure and regularity of the global attractor of a reacction-diffusion
equation with non-smooth nonlinear term, \emph{Discrete Continuous Dynamical
Systems,} \textbf{32} (2014), 4155-4182.

\bibitem {kapustyankasyanovvalero2015}O.V. Kapustyan, P.O. Kasyanov, J.
Valero, Structure of the global attractor for weak solutions of a
reaction-diffusion equation, \emph{Applied Mathematics \& Information
Sciences,} \textbf{9} (2015), 2257-2264.

\bibitem {kapustyanpankov}O. V. Kapustyan, V. Pankov and J. Valero, On global
attractors of multivalued semiflows generated by the 3D B\'enard system,
\emph{Set-Valued and Variational Analysis}, \textbf{20} (2012), 445-465.

\bibitem {li}D. Li, Morse decompositions for general dynamical systems and
differential inclusions with applications to control systems, \emph{SIAM
Journal on Control and Optimization}, \textbf{46} (2007), 35-60.

\bibitem {lions}J. L. Lions, {\itshape Quelques M\'ethodes de R\'esolution des
Probl\`emes aux Limites non Lin\'eaires}, Gauthier-Villar, Paris, 1969.

\bibitem {simone}S. Mazzini, \emph{Atratores para o problema de
Chafee-Infante}, PhD-thesis, Universidad de S\~ao Paulo, 1997.

\bibitem {melnikvalero}V. S. Melnik and J. Valero, On attractors of
multi-valued semi-flows and differential inclusions, \emph{Set-Valued
Analysis}, \textbf{6} (1998), 83-111.

\bibitem {robinson}J. C. Robinson, {\ \itshape Infinite-Dimensional Dynamical
Systems: An Introduction to Dissipative Parabilic PDEs and the Theory of
Global Attractors}, Cambridge University Press, Cambridge, UK, 2001.

\bibitem {sellyou}G. R. Sell and Y. You, {\itshape Dynamics of evolutionary
equations}, Springer, 2002.

\bibitem {temmam}R. Temam, {\itshape Navier-Stokes Equations}, North-Holland,
Amsterdam-New York, 1977.

\bibitem {tolstonogov}A. Tolstonogov, On solutions of evolution inclusions I,
\emph{Siberian Math. J.}, \textbf{33} (1992), 500-511.

\bibitem {valero2001}J. Valero, Attractors of parabolic equations without
uniqueness, \emph{J. Dynamics Differential Equations}, \textbf{13} (2001), 711-744.

\bibitem {valero2005}J. Valero, On the Kneser property for some parabolic
problems, \emph{Topology Appl.,} \textbf{153} (2005), 975-989.

\bibitem {valerokapustyan}J. Valero and A. V. Kapustyan, On the connectedness
and symptotic behaviour of solutions of reaction-diffusion systems, \emph{J.
Math. Anal. Appl.}, \textbf{323} (2006), 614-633.

\bibitem {daleroberts}A. Wayne and D. Varberg, {\itshape Convex functions},
Academic Press, Elsevier, 1973.

\bibitem {yosida}K. Yosida, {\itshape Functinoal Analysis}, Springer-Verlag,
Berlin, 1965.
\end{thebibliography}
\end{document}